\documentclass[smallextended,referee,envcountsect,]{svjour3}
\smartqed
\usepackage{graphicx}
\usepackage{color}
\usepackage{amsmath}
\usepackage{bm}
\usepackage{amssymb}
\usepackage{enumitem}
\usepackage{hyperref}
\usepackage{cleveref}
\usepackage{verbatim}
\setlist[enumerate]{label=(\roman*)}
\usepackage{booktabs}

\DeclareMathOperator*{\argmax}{arg\,max}
\journalname{JOTA}

\begin{document}

\title{On Completely Positive and Semidefinite Programs for Moment Bounds on Maximum of Random Quadratics}

\subtitle{}

\author{Buse Unluturk \and Divya Padmanabhan \and Karthik Natarajan}

\institute{Buse Unluturk and Karthik Natarajan \at
Engineering Systems Design, Singapore University of Technology and Design, 8 Somapah Road, Singapore 487372\\
buse\_unluturk@sutd.edu.sg,
karthik\_natarajan@sutd.edu.sg
\and
Divya Padmanabhan \at
School of Mathematics and Computer Science,
Indian Institute of Technology (IIT) Goa, Farmagudi, Ponda, Goa 403401, India.\\
divya@iitgoa.ac.in
}
%\institute{Tam\'as Terlaky \at
%             Lehigh University \\
%              Bethlehem, PA 18015, USA\\
%              terlaky@Lehigh.edu
%              Virtual Author,  Corresponding author  \at
%              University of Dreamland \\
%              Dreamland\\
%              author2@example.com
%}

\date{Received: date / Accepted: date}
%The correct dates will be entered by the editor.

\maketitle

\begin{abstract}
Completely positive programs are known to provide an exact representation to several NP-hard optimization problems. In this work, we look at the maximum of random quadratic functions and compute the maximum expected value over ambiguity sets defined with moment and support information. We propose exact completely positive programs and semidefinite programs when the random variables are Bernoulli, nonnegative and unrestricted in support. Our experimental study demonstrates the utility of these formulations.
%Each submission must include an abstract, usually between 50 to 200 words. The abstract should give a clear and concise idea of the  content of the paper,  without the use of mathematical formulas, equations, and acronyms.
\end{abstract}
%\subclass{49J53 \and  49K99 \and more}

%All acknowledgements should be placed in the back of the paper after Conclusions..

\section{Introduction}
\label{sec:introduction}
In this work, we study the maximum of random quadratic functions of the form,
$
  f(\bm{\tilde{\xi}} )=  \max_{k \in [K]} \bm{\tilde{\xi}}^T\mathbf{A}_k \bm{\tilde{\xi}}+\bm{\tilde{\xi}}^T\bm{b}_k +c_k
$
where $[K] = \{ 1, \ldots, K \}$ and $\bm{\tilde{\xi}}$ is a random vector of dimension $n$. For any  $k \in [K]$, $\mathbf{A}_k$ is a real valued symmetric matrix of size $n\times n$, $\bm{b}_k$ is a vector in $\mathbb{R}^n$ and $c_k$ is  a scalar. 
%and hence $f(\bm{\tilde{\xi}})$ is a random piecewise quadratic function. 
%The matrices $\mathbf{A}_k$ are not restricted to  have any particular structure. 
The probability distribution of the random vector $\bm{\tilde{\xi}}$ is not precisely known and assumed to lie in an ambiguity set denoted by $\mathcal{P}$. We are interested in computing the maximum expected value among all feasible distributions in the ambiguity set: 
\begin{align}
\label{sup_problem}
\tag{P}
    \sup_{\mathbb{P} \in \mathcal{P}} \mathbb{E}_{\mathbb{P}}\left[\max_{k \in [K]} \left(\bm{\tilde{\xi}}^T\mathbf{A}_k \bm{\tilde{\xi}}+\bm{\tilde{\xi}}^T\bm{b}_k +c_k\right)\right].
\end{align}

With focus on problem \eqref{sup_problem}, we restrict our attention to the case where the ambiguity set $\mathcal{P}$ consists of all distributions of the random vector with given means, variances and lower bounds on covariances. This moment ambiguity set has been recently proposed in \cite{natarajan2025submodulardro}.
%and describes random vectors with known first and second moments and fixed lower bounds on cross moments (equivalenly fixed means, variances and lower bounds on covariances).  
We consider three cases where the support of the random vector is either  (a) $\{0,1\}^n$, (b) $\mathbb{R}_n^+$ and (c) $\mathbb{R}_n$. As shown in \cite{natarajan2025submodulardro} checking for the non-emptiness of the ambiguity set can be done by solving a linear program in case (a) and by solving semidefinite programs in cases (b) and (c). Here we develop moment bounds for this ambiguity set for general $\mathbf{A}_k$ matrices while in  \cite{natarajan2025submodulardro}, the results were developed for restricted $\mathbf{A}_k$ matrices with nonnegative off diagonal entries.
\\

\noindent\textbf{Motivation} 
The problem is specially of relevance in distributionally robust optimization when there is a decision $\bm{x}$ to be made and the matrices $\mathbf{A}_k$, the vectors $\bm{b}_k$ and the scalars $c_k$ are all functions of $\bm{x}$. When $\bm{\xi}$ is deterministic and known, one seeks to find $\bm{x}$ so that,
\begin{align*}
    \inf_{\bm{x} \in \mathcal{X}} f(\bm{\xi}, \bm{x})= \inf_{\bm{x} \in \mathcal{X}}\max_{k\in[K]}\left(\bm{{\xi}}^T\mathbf{A}_k(\bm{x})\bm{{\xi}}+\bm{{\xi}}^T\bm{b}_k(\bm{x})+c_k\right)
\end{align*} However, in the presence of uncertainty in $\bm{\tilde{\xi}}$, one possible approach is to solve the following distributionally robust optimization problem,
\begin{align*}
    \inf_{x\in\mathcal{X}} \sup_{\mathbb{P}\in\mathcal{P}} \text{E}_{\mathbb{P}}\left[\max_{k\in[K]}\left(\bm{\tilde{\xi}}^T\mathbf{A}_k(\bm{x})\bm{\tilde{\xi}}+\bm{\tilde{\xi}}^T\bm{b}_k(\bm{x})+c_k\right)\right]
\end{align*}
%When $K=1$, $f(\bm{\tilde{\xi}}, \bm{x})$ is a random quadratic function. Quadratic functions show up in several avenues. 
Such problems show up for example, in regression tasks in machine learning where we have a noisy training dataset and the square loss (a quadratic function) is to be minimized. The training dataset consists of noisy features often, in which one desires a decision $\bm{x}$ to be learnt that is robust to the values of the features given in the training  dataset where $\bm{\tilde{\xi}}$ is a function of the training dataset (which is random).

Our key contributions are the following:
\begin{enumerate}
\item We show in  \Cref{lemma:nphard} that the problem \eqref{sup_problem} is NP-hard for the ambiguity set we consider. Specifically the case of Bernoulli support and $K=1$ is used for the hardness proof.
    \item We then propose exact conic programming formulations for solving \eqref{sup_problem} over the moment ambiguity set. For the case of Bernoulli support and $\mathbb{R}_n^+$, the proposed formulations are completely positive programs while for support in $\mathbb{R}_n$, we propose a semidefinite program. We prove the tightness of these formulations in \Cref{thm:bernoulli} and \Cref{thm:R_n_+_and_R_n}.
    \item Our numerical results study the optimality gap of the doubly nonnegative relaxation of the completely positive program for the Bernoulli support. The exact bound computed via an exponential sized linear program becomes very expensive as the number of random variables grow. On the other hand, the doubly nonnegative relaxation though not exact provides a computationally tractable procedure to obtain useful bounds.  
\end{enumerate}

The rest of the paper is organized as follows. In \Cref{sec:litReview}, we discuss related work and also introduce the reader to the cone of completely positive matrices and related cones. In \Cref{sec:conicProgrammingFormulation}, we propose our conic programming formulations for \eqref{sup_problem} while tightness proofs are provided in \Cref{sec:tightnessProofs}. Finally we describe our numerical results in \Cref{sec:numericals} and conclude in \Cref{sec:conclusion}. 
\\

\noindent \textbf{Notations:} We denote scalar values as   $\lambda$, vectors with bold letters as $\bm{\lambda}$, matrices with bold capital letters as $\mathbf{\Lambda}$ and random variables and random vector with the tilde notation as $\tilde{\xi}$ and $\tilde{\bm{\xi}}$. We use $[n]$ to denote the set $\{1, \ldots, n\}$ and $\binom{[n]}{2}$ to denote the set $\{(1,2),(1,3),\ldots,(n-1,n)\}$ of size $\binom{n}{2}$. We use $\mathbf{1}$ to denote a vector of all ones. $\mathbb{R}_n$ denotes the set of $n$ dimensional real vectors while $\mathbb{R}_n^+$ denotes the set of nonnegative vectors of size $n$. $\mathbf{A} \succeq 0$ denotes a positive semidefinite matrix $\mathbf{A}$ and $\mathbf{A} \geq 0$ denotes a nonnegative matrix $\mathbf{A}$. The trace of a matrix $\mathbf{A}$ is denoted as $trace(\mathbf{A})$. The inner product  between two matrices of the same dimension is denoted as $\mathbf{A}\bullet\mathbf{\Lambda}=trace(\mathbf{A}^T\mathbf{\Lambda})$. The diagonal entries of a square matrix $\mathbf{A}$ are denoted by $\mbox{diag}(\mathbf{A})$. $\mathbb{R}_{n\times n}^+$ is the set of $n\times n$ nonnegative matrices, and $\mathbb{S}_n$ is the set of $n\times n$ symmetric matrices. $\mathbb{I}(\text{E})$ is an indicator function that takes the value 1 when the event E occurs, 0 otherwise.

\section{Literature Review}
\label{sec:litReview}
We first begin our literature review by discussing the completely positive cone and its application in (distributionally) robust optimization. 
\subsection{Completely Positive Cone and Approximation Hierarchies}
A completely positive matrix $\bm{X}$ of size $n \times n$ (denoted $\bm{X} \succeq_{cp} 0$) is a symmetric matrix which admits a decomposition of the form  $\bm{X}  = \bm{Z}\bm{ Z}^\top \text{ for some } \bm{Z} \in \mathbb{R}_{n \times k}^+ \text{ where } k$ is some positive integer. An equivalent characterization is the existence of $k$ nonnegative vectors $\left\lbrace \bm{v}_1, \ldots, \bm{v}_k \right \rbrace$ of size $n$ such that $\bm{X}  = \sum_{i=1}^k \bm{v}_i \bm{v}_i^\top$ for some value of $k$.
The completely positive cone is a set of all completely positive matrices defined as
$
    C_n = \left \lbrace \bm{X} \in \mathbb{S}_{n}: \bm{X}  \succeq_{cp} 0 \right \rbrace
$. It is easy to see that $C_n$ is a cone. The dual cone to the completely positive cone is the set of copositive matrices defined as:
\begin{align*}
    C^*_n = \left \lbrace \bm{X}  \in \mathbb{S}_{n}: \bm{v}^{\top} \bm{X}  \bm{v} \geq 0, \ \forall \bm{v} \in \mathbb{R}_n^+ \right \rbrace.
\end{align*} Testing for membership in $C_n$ is NP-hard while testing for membership in $C^*_n$ is co-NP-hard \cite{dickinson2014computational}. A related cone is the cone of  doubly nonnegative matrices defined as, 
\begin{align*}
    D_n = \left \lbrace \bm{X}  \in \mathbb{S}_{n}: \bm{X}  \succeq 0, \bm{X} \geq 0\right \rbrace.
\end{align*}
For $n \leq 4$, it is well known that $C_n = D_n$, however for $n > 4, $ $C_n \subset D_n$.

A completely positive program  is a linear optimization problem over the cone $C_n$ given by \begin{align*}
    \min (\max)\, \left\{\bm{C} \bullet \bm{X}  \ | \ \bm{G} _i\bullet \bm{X} = g_i, \ \forall i \in [m],\bm{X}  \in \text{C}_n\right\}.
    \end{align*}
    Since testing membership in $C_n$ is itself hard, solving a completely positive program is also NP-hard. The dual of the completely positive program  is referred to as  copositive program and involves an optimization over $C^*_n$,  which is also a hard problem. However despite the hardness, tractable relaxations for the completely positive program exist through semidefinite programming approaches by replacing the constraint $\bm{X} \in C_n$ with $\bm{X} \in D_n$. Such relaxations are referred to as doubly nonnegative relaxations (DNN).
For more information and properties of $C_n$, $D_n$ and completely positive/copositive optimization, the reader may refer to the excellent textbook on completely positive matrices \cite{Berman2003}.
\subsection{Applications in (Distributionally) Robust Optimization}
Several works have looked at the problem of maximizing the expected value of piecewise linear functions under various assumptions on the set $\mathcal{P}$ (the objective then is a special case of our setting in \eqref{sup_problem} with $\bm{A}_k = 0$ for all $k$). For support in $\mathbb{R}_n$ and when  first two moments are enforced exactly,  the problem \eqref{sup_problem} can be solved exactly using an SDP involving $K$ matrices each of size $n+1$ \cite{Mishra2012,Natarajan2017}. 

A closely related setting is when the objective function $f(\bm{\tilde{\xi}}) = \max_{\bm{z} \in \mathcal{Z}} \bm{\tilde{\xi}}^{\top} \bm{z}$ where $\mathcal{Z}$ is a polytope. One can represent  this function $f$ in the form given in \Cref{sec:introduction} by setting $K$ as the number of extreme points of $\mathcal{Z}$ and by representing the extreme points using the vectors $\bm{b}_k$. \cite{Natarajan2011} proposed formulations for the case of uncertain mixed integer linear programs represented via the polytope $\mathcal{Z}$. The support considered for the random vector therein was $\mathbb{R}_n^+$ and the first two moments of $\bm{\xi}$ were assumed to be known perfectly. For this case an exact completely positive program involving a matrix of size $(2n+1) \times (2n+1)$ was proposed. Their approach built on completely positive formulations for non-convex QPs which has in turn received much attention in the last three decades \cite{bomze2002solving,Burer09,bomze2017fresh}.
A key quantity which shows up in the approach of \cite{Natarajan2011} is the random vector $\bm{q(\tilde{\xi})} = [1, \bm{\tilde{\xi}}, \bm{\tilde{z(\xi)}}]$, where $\bm{\tilde{z(\xi)}}$ is the maximizer of $\bm{\tilde{\xi}}^\top \bm{z}$ and hence is random. The reformulations in \cite{Natarajan2011} are constructed by characterizing the first and second moments of $\bm{q(\tilde{\xi})}$ in terms of known information on the moments of $\bm{\tilde{\xi}}$. Similar approach has been used in \cite{padmanabhan2021exploiting} to derive SDPs involving smaller sized matrices when block diagonal submatrices of $\mathbb{E}[\tilde{\bm{\xi}} \tilde{\bm{\xi}}^\top]$ are known exactly. These formulations have been used in a wide variety of applications including appointment scheduling, newsvendor problem, discrete choice, process flexibility and even gaming  \cite{kong2013scheduling,Natarajan2009,padmanabhan2021exploiting,liu2024picking,liu2023picking,Yan2018}. While we enforce lower bounds on the cross terms in the covariance matrix, the set in \cite{Delage2010} enforces a matrix lower bound in semidefinite order. This boils down to enforcing the matrix $\text{E}_\mathbb{P}[\tilde{\bm{\xi}}\tilde{\bm{\xi}}^{\top}]-\bm{\Sigma}$ to be positive semidefinite. A detailed treatment of various formulations arising in the context of moment based ambiguities can be seen in the book \cite{karthikBook}.
On a related note \cite{Hanasusanto2018}
propose copositive reformulations for two-stage distributionally robust linear programs where the ambiguity set $\mathcal{P}$ is a  Wasserstein Ball based set. Here $\mathcal{P}$  consists of all distributions within a specified Wasserstein distance from an empirical distribution constructed using available samples.

While all the aforementioned references look at linear or piecewise linear functions in the objective, the case of a single quadratic objective has been looked at in \cite{Hanasusanto2024} in the larger context of a distributionally robust two stage stochastic program. One of the subproblems looked at in the paper is $\sup \{\mathbb{E}[\bm{\xi}^\top A \bm{\xi}] \ | \ \mathbb{P} \in \mathcal{P}\}$. The ambiguity set $\mathcal{P}$ assumes availability of samples and the  Frobenius norm of the difference between the true second moment and empirical second moment is assumed to be bounded by a parameter $\gamma$. Under this setup, a completely positive optimization formulation is proposed. Our framework on the other hand assumes known lower bounds on the second moment matrix and also considers maximum of quadratic functions in the objective.  

Copositive programming approaches have also been used in the context of robust two stage stochastic programs and adjustable robust optimization \cite{xu2025improved,xu2018copositive,bomze2023adjustable}. In robust optimization, instead of distributional uncertainty, the uncertain parameter $\bm{\xi}$ is assumed to take values from a set $\Xi$ and one solves a problem of the form 
$\inf_{\bm{x} \in \mathcal{X}} \sup_{\bm{\xi} \in \Xi} f(\bm{\xi}, \bm{x})$. In general completely positive and copositive formulations are useful as reformulations for NP-hard problems since their approximations via DNN relaxation are implementable.

%It is well known that solving a non-convex QP is NP-hard in general.  The theory of completely positive programming provides alternate representations for quadratic programs (not necessarily convex).  Several papers have been written  towards solving such QPs including  These formulations solve a linear optimization problem over the cone CPP$_n$ of completely positive matrices \cite{Berman2003}, 

%In the last two decades, the above theory on QPs has been adapted to solve uncertain mixed integer linear optimization problems under various ambiguities. When the first two moments of coefficient vector are known, one can construct a set of distributions satisfying these moments. The underlying separation problem for these uncertain problems involve QPs and hence the theory of QPs is of relevance. In particular \cite{Natarajan2011} assume support in $\mathbb{R}^n_+$ and propose a completely positive program over the cone CPP$_{2n+1}$. On the other hand when the support of the random variables lie in $\mathbb{R}^n$ a tractable semidefinite program provides a tight formulation \cite{Natarajan2017} and one does not need to go into completely positive cone.\\

\section{Conic Programming Reformulations}
\label{sec:conicProgrammingFormulation}
In this section, we propose a general conic program to solve \eqref{sup_problem}. The ambiguity set is defined as follows:
\begin{align}
    \mathcal{P}=\left\{\mathbb{P}\in\mathcal{P}\left(\Xi\right)\ \big| \ \mathbb{E}_\mathbb{P}[\bm{\tilde{\xi}}] = \bm{\mu},  \mathbb{E}_\mathbb{P}[\mbox{diag}(\bm{\tilde{\xi}}\bm{\tilde{\xi}}^T)] = \mbox{diag}(\mathbf{\Sigma}),      \mathbb{E}_\mathbb{P}[\bm{\tilde{\xi}}\bm{\tilde{\xi}}^T]\geq \mathbf{\Sigma}
 \right\},
 \label{ambiguity_set}
\end{align}
where $\mu_i$ is the given mean of ${\tilde{\xi}}_i$, $\Sigma_{ii}$ is the given second moment of $\tilde{\xi}_i$ and $\Sigma_{ij}$ is the given lower bound on the cross moment of $\tilde{\xi}_i\tilde{\xi}_j$.

Before we propose our conic program, we first
show that the problem is generally NP-hard even for the case of $\Xi = \{0,1\}^n$. 
\begin{lemma}
\label{lemma:nphard}
The problem \eqref{sup_problem} is NP-hard even for the case where $\Xi = \{0,1\}^n$.
\end{lemma}
{\it Proof.}
We will prove the result for the case of $K=1$. Then let $\mathbf{A} = \mathbf{A}_1, \mathbf{b} = \mathbf{b}_1$ and $c = c_1$. Denote by $\theta(\bm{\xi}) = \mathbb{P}(\bm{\tilde{\xi}} = \bm{\xi})$. Then the problem \eqref{sup_problem} can be written as a linear program over $2^n$ decision variables as follows,
\begin{align}
\max_{\theta} &\, \sum_{\bm{\xi} \in \{0,1\}^n}\theta(\bm{\xi}) (\bm{\xi}^T\mathbf{A}\bm{\xi}+\bm{\xi}^T\bm{b} +c) \label{exp_size_lp_k1} \\
\text{s.t} &\, \sum_{\bm{\xi} \in \{0,1\}^n}\theta(\bm{\xi}) = 1, \nonumber\\
&\, \sum_{\bm{\xi} \in \{0,1\}^n}\theta(\bm{\xi}) \xi_i = \mu_i, \ \forall i \in [n],  \nonumber \\
%&\, \sum_{\bm{\xi} \in \{0,1\}^n}\theta(\bm{\xi}) \xi_i^2 = \Sigma_{ii}, \ \forall i \in [n],  \nonumber \\
&\, \sum_{\bm{\xi} \in \{0,1\}^n}\theta(\bm{\xi}) \xi_i \xi_j \geq \Sigma_{ij}, \ \forall (i,j) \in \binom{[n]}{2},\\
&\, \theta(\bm{\xi}) \geq 0, \ \forall \bm{\xi} \in \{0,1\}^n.
\nonumber
\end{align}
Here we omit the constraint $\text{E}_\mathbb{P}[\mbox{diag}(\bm{\tilde{\xi}}\bm{\tilde{\xi}}^T)] = \mbox{diag}(\mathbf{\Sigma})$ since for Bernoulli random variables, $\tilde{\xi}_i^2 = \tilde{\xi}_i \forall i \in [n]$ and we must have $\mbox{diag}(\mathbf{\Sigma}) = \bm{\mu}$. The dual of \eqref{exp_size_lp_k1} can be written as,
\begin{align*}
    \min &\, \rho + \bm{\alpha}^{\top} \bm{\mu} + \bm{\Gamma} \bullet \Sigma \\
    &\, \rho + \bm{\alpha}^{\top} \bm{\xi} + \bm{\xi}^{\top}  \bm{\Gamma} \bm{\xi} \geq c + \bm{b}^{\top} \bm{\xi} + \bm{\xi}^{\top}  \bm{A}  \bm{\xi}, \  \forall \bm{\xi} \in \{0,1\}^n\\
    & \,\mbox{diag}(\bm{\Gamma}) = \bm{0}\\
    &\,  \bm{\Gamma} \leq 0.
\end{align*}
%In the dual note that $\rho \in \mathbb{R}, \bm{\alpha} \in \mathbb{R}_n$, $\Gamma(i,j) \leq 0 \,\forall i,j \in \binom{[n]}{2}$ and $\Gamma(i,i) = 0 \,\forall i \in [n]$. 
The separation problem thus boils down to,
\begin{align}
\label{separation_problem}
    \min_{\bm{\xi} \in \{0,1\}^n} (\rho -c ) + (\bm{\alpha} - \bm{b})^{\top}\bm{\mu} + \bm{\xi}^{\top} (\bm{\Gamma-A})\bm{\xi} 
\end{align}
with $\rho, c, \bm{\alpha}, \bm{b}, \bm{\Gamma}, \bm{A}$ fixed. Since the matrix $\bm{A}$ does not have any particular structure, the separation problem is an instance of a QUBO (quadratic unconstrained binary optimization) problem which is known to be NP-hard \cite{gareyJohnson}. Hence hardness of $\eqref{sup_problem}$ is proved by equivalence of separation and optimization \cite{GLS93}.
\qed
We note that while the problem we consider in this paper is generally NP-hard, the case where the off-diagonal entries of each $\mathbf{A}_k$ are nonnegative is polynomial time solvable as proved in \cite{natarajan2025submodulardro}. The difference it would make in the proof approach in \Cref{lemma:nphard} (for $K=1)$ is that due to non-negativity of the cross terms in $\mathbf{A}$, cross terms of $\bm{\Gamma - A}$ are negative, and \eqref{separation_problem} would boil down to a polynomial time solvable case of QUBO \cite{bqp_poly}. Another related setting is an ambiguity set where univariate marginal distributions on $\bm{\tilde{\xi}}$ are known. For the case of Bernoulli random variables, this would boil down to an ambiguity set  similar to \Cref{ambiguity_set} where only the first moment $\bm{\mu}$ is known and information on the second moment is unknown. This setting would also be NP-hard due to a very similar proof argument as in \Cref{lemma:nphard}. For $K=1$, the separation problem would boil down to  
\begin{align}
\label{separation_problem_marginal}
    \min_{\bm{\xi} \in \{0,1\}^n} (\rho -c ) + (\bm{\alpha} - \bm{b})^{\top}\bm{\mu} -  \bm{\xi}^{\top} \bm{A}  \bm{\xi}
\end{align}
which is also NP-hard since $\bm{A}$ is not structured in any way.

We will now provide formulations for the following three cases,
(a) $\Xi = \{ 0, 1\}^n$, 
(b) $\Xi = \mathbb{R}_n^+$ and
(c) $\Xi= \mathbb{R}_n $. 
A few of the decision variables appear commonly across all the three support sets  and we define them first, along with the constraints to be enforced on them. For any $\bm{\tilde{\xi}}$, suppose $k^*(\bm{\tilde{\xi}}) \in  \argmax_{k \in [K]} \bm{\tilde{\xi}}^T\mathbf{A}_k \bm{\tilde{\xi}}+\bm{\tilde{\xi}}^T\bm{b}_k +c_k$. If there are multiple optimal indices at any realization $\bm{\xi}$, we let $k^*(\bm{\xi})$ be a measurable selection on the set of optimal indices. Since $\bm{\tilde{\xi}}$ is random, $k^*(\bm{\tilde{\xi}})$ is also random. Over the distribution governing $\bm{\tilde{\xi}}$, define,
\begin{align*}
    \lambda_{0,k}&=\mathbb{P}(k^*(\bm{\tilde{\xi}})=k), \quad \forall k \in [K],\\
    \bm{\lambda}_{1,k}&=\mathbb{E}\left[\bm{\tilde{\xi}}\mathbb{I}\left(k^*(\bm{\tilde{\xi}})=k\right)\right], \quad \forall k \in [K],\\
\mathbf{\Lambda}_{1,k}&=\mathbb{E}\left[\bm{\tilde{\xi}}\bm{\tilde{\xi}}^T\mathbb{I}\left(k^*(\bm{\tilde{\xi}})=k\right)\right], \quad \forall k \in [K].
\end{align*}
%where $\mathbb{P}\left(k^*(\bm{\tilde{\xi}})=k\right)$ is the probability of the optimal index $k^*(\bm{\tilde{\xi}})$ being $k$. 
Then, the objective function can be written  as
\begin{align*}
    \mathbb{E}_\mathbb{P} &\left[\max_{k\in[K]}\left(\bm{\tilde{\xi}}^T\mathbf{A}_k\bm{\tilde{\xi}}+\bm{\tilde{\xi}}^T\bm{b}_k+c_k\right)\right]=\mathbb{E}_\mathbb{P}\left[\sum_{k\in[K]}\left(\bm{\tilde{\xi}}^T\mathbf{A}_k\bm{\tilde{\xi}}+\bm{\tilde{\xi}}^T\bm{b}_k+c_k\right)\mathbb{I}(k^*(\bm{\tilde{\xi}})=k)\right]\\
    &=\sum_{k\in[K]} \mathbb{E}_\mathbb{P}\left[\left(\bm{\tilde{\xi}}^T\mathbf{A}_k\bm{\tilde{\xi}}+\bm{\tilde{\xi}}^T\bm{b}_k+c_k\right)\mathbb{I}(k^*(\bm{\tilde{\xi}})=k)\right]\\
    %&=\sum_{k\in[K]}\left[\sum_{i=1}^n\sum_{j=1}^n A_k(i,j)\mathbb{E}_\mathbb{P}[\tilde{\xi}_i\tilde{\xi}_j\mathbb{I}(k^*(\bm{\tilde{\xi}})=k))]
    %+\sum_{i=1}^n b_k(i)\mathbb{E}_\mathbb{P}[\tilde{\xi}_i\mathbb{I}(k^*(\bm{\tilde{\xi}})=k)]
    %+c_kP(k^*(\bm{\tilde{\xi}})=k)\right]\\
&=\sum_{k\in[K]}\left[\mathbf{\Lambda}_{1,k}\bullet\mathbf{A}_k+\bm{\lambda}_{1,k}^T\bm{b}_k+\lambda_{0,k}c_k\right]
\end{align*}
The mean of the random variables in $\mathcal{P}$ can be expressed as:
\begin{align*}
    \begin{array}{lll}
       \displaystyle  \mathbb{E}[\bm{\tilde{\xi}}] & \displaystyle =\sum_{k\in[K]} \mathbb{E}[\bm{\tilde{\xi}}\mathbb{I}(k^*(\bm{\tilde{\xi}})=k)]= \sum_{k\in[K]}\bm{\lambda}_{1,k} = \bm{\mu}.
    \end{array} 
    %\quad \Rightarrow \quad  \sum_{k\in[K]}\bm{\lambda}_{1,k} = {\bm{\mu}}.
\end{align*}
The second moments of the random variables in $\mathcal{P}$ can be expressed as:
\begin{align*}
    \begin{array}{lll}
    \displaystyle    \mathbb{E}[\bm{\tilde{\xi}}\bm{\tilde{\xi}}^T] & \displaystyle=\sum_{k\in[K]}\mathbb{E}\left[\bm{\tilde{\xi}}\bm{\tilde{\xi}}^T\mathbb{I}(k^*(\bm{\tilde{\xi}})=k)\right]
    =\sum_{k\in[K]}\mathbf{\Lambda}_{1,k}.
        \end{array} 
\end{align*}
This gives the constraints
\begin{align*}
    \begin{array}{lll}
      &  \displaystyle \sum_{k\in[K]}\mbox{diag}(\mathbf{\Lambda}_{1,k})=\mbox{diag}(\mathbf{\Sigma}),\\ 
    &  \displaystyle \sum_{k\in[K]}\mathbf{\Lambda}_{1,k}\geq\mathbf{\Sigma}. 
       \end{array} 
\end{align*}
Since the probabilities over all values of $k^*(\bm{\tilde{\xi}})$ must sum to 1, we must have
$
    \sum_{k\in[K]}\lambda_{0,k}=1.$

We are now ready to propose our conic  formulation. For the cases of $\Xi = \{0,1\}^n, \Xi = \mathbb{R}_n^+$ and $\Xi = \mathbb{R}_n$,  we propose the following  conic programming formulation to solve \eqref{sup_problem}.
\begin{align}
\label{CPP}
    \tag{CPP} \quad \quad \max& \sum_{k\in[K]}\left[\mathbf{\Lambda}_{1,k}\bullet\mathbf{A}_k+\bm{\lambda}_{1,k}^T\bm{b}_k+\lambda_{0,k}c_k\right] \\
    \text{s. t.}&\label{C1}\tag{C1} \quad \sum_{k\in[K]}\bm{\lambda}_{1,k} = \bm{\mu}, \\
    &\label{C2}\tag{C2} \quad \sum_{k\in[K]}\mbox{diag}(\mathbf{\Lambda}_{1,k})= \mbox{diag}(\mathbf{\Sigma}) \\
    &\label{C3}\tag{C3} \quad \sum_{k\in[K]}\mathbf{\Lambda}_{1,k}\geq \mathbf{\Sigma}, \\
    &\label{C4}\tag{C4} \quad \sum_{k\in[K]}\lambda_{0,k}=1, \\
    & \left\{\left(\lambda_{0,k}, \bm{\lambda}_{1,k}, \mathbf{\Lambda}_{1,k}\right)\right\}_{k\in[K]}\in\mathcal{F} \nonumber
\end{align}
where the specification of the set $\mathcal{F}$  varies according to the support set under consideration. For all the cases, $\mathcal{F}$ turns out to be a convex cone. In \Cref{subsec:01support,subsec:rnplus_support,subsec:rn_support} we describe this set $\mathcal{F}$ for each of the support sets we consider.
\subsection{Support $\Xi = \{0,1\}^n$}
\label{subsec:01support}
For the support $\{0,1\}^n$, we have $\Sigma_{ii} = \mu_i$ for each $i \in [n]$. In addition, we require additional decision variables. Define
\begin{align*}
    \bm{\lambda}_{2,k}=&\mathbb{E}[(\mathbf{1}-\bm{\tilde{\xi}})\mathbb{I}(k^*(\bm{\tilde{\xi}})=k)], \quad \forall k \in [K]\\
    \mathbf{\Lambda}_{2,k}=&\mathbb{E}[(\mathbf{1}-\bm{\tilde{\xi}})\bm{\tilde{\xi}}^T\mathbb{I}(k^*(\bm{\tilde{\xi}})=k)], \quad \forall k \in [K]\\
    \mathbf{\Lambda}_{3,k}=&\mathbb{E}[(\mathbf{1}-\bm{\tilde{\xi}})(\mathbf{1}-\bm{\tilde{\xi}})^T\mathbb{I}(k^*(\bm{\tilde{\xi}})=k)],\quad \forall k \in [K].
\end{align*}
Then,  $\mathcal{F}$ is the set,
\begin{align*}
    \mathcal{F} = \{(\lambda_{0,k}, \bm{\lambda}_{1,k}, \mathbf{\Lambda}_{1,k}) \text{ } \forall k\in[K]: \exists (\bm{\lambda}_{1,k}, \mathbf{\Lambda}_{2,k},
        \mathbf{\Lambda}_{3,k}),
        \forall k\in[K] &\\
        \text{ satisfying constraints (\ref{C5},\ref{C6},\ref{C7},\ref{C8})  below.} \}
\end{align*}

\begin{align}
         \label{C5}\tag{C5}\quad & \sum_{k\in{K}}\left(\bm{\lambda}_{1,k}+\bm{\lambda}_{2,k}\right)=\mathbf{1}, \\
        \quad & \sum_{k\in[K]}\left(\mbox{diag}(\mathbf{\Lambda}_{1,k})+2\mbox{diag}(\mathbf{\Lambda}_{2,k})+\mbox{diag}(\mathbf{\Lambda}_{3,k})\right)= \mathbf{1},\label{C6}\tag{C6}\\
        \label{C7}\tag{C7} \quad & \sum_{k\in[K]}\bm{\lambda}_{1,k}=\sum_{k\in[K]}\mbox{diag}(\mathbf{\Lambda}_{1,k}),\\
        \label{C8}\tag{C8} \quad &\begin{bmatrix}
            \lambda_{0,k} & \bm{\lambda}_{1,k}^T & \bm{\lambda}_{2,k}^T \\
            \bm{\lambda}_{1,k} & \mathbf{\Lambda}_{1,k} & \mathbf{\Lambda}_{2,k}^T \\
            \bm{\lambda}_{2,k} & \mathbf{\Lambda}_{2,k} & \mathbf{\Lambda}_{3,k}
        \end{bmatrix}\in \text{C}_{2n+1}\, \forall k \in [K]
     %   \label{C9}\tag{C9} \quad & \bm{\lambda}_{2,k}\in\mathbb{R}_n^+, \mathbf{\Lambda}_{2,k}\in\mathbb{R}_{n\times n}^+, \mathbf{\Lambda}_{3,k}\in\mathbb{S}_n\cap\mathbb{R}_{n\times n}^+ \quad \forall k \in [K]
    \end{align}
    We observe that while matrices $\mathbf{\Lambda}_{1,k}$ and $\mathbf{\Lambda}_{3,k}$ are symmetric square matrices, $\mathbf{\Lambda}_{2,k}$ is square but need not be symmetric.
\iffalse
\begin{align}
    \left\{
    \begin{matrix}
        \left(\lambda_{0,k}, \bm{\lambda}_{1,k},
        \mathbf{\Lambda}_{1,k}\right),\\
        \forall k\in[K]
    \end{matrix}: 
    \begin{aligned}
         \quad & \sum_{k\in{K}}\left(\lambda_{1,k}(i)+\lambda_{1,k}(i)\right)=1 \quad \forall i\in[n], \\
        \quad & \sum_{k\in[K]}\left(\Lambda_{1,k}(i,i)+2\Lambda_{2,k}(i,i)+\Lambda_{3,k}(i,i)\right)=1 \quad \forall i \in[n],\\
        \text{(C6)} \quad & \sum_{k\in[K]}\lambda_{1,k}(i)=\sum_{k\in[K]}\Lambda_{1,k}(i,i)\quad \forall i \in[n],\\
        \text{(C7)} \quad &\begin{bmatrix}
            \lambda_{0,k} & \bm{\lambda}_{1,k}^T & \bm{\lambda}_{2,k}^T \\
            \bm{\lambda}_{1,k} & \mathbf{\Lambda}_{1,k} & \mathbf{\Lambda}_{2,k}^T \\
            \bm{\lambda}_{2,k} & \mathbf{\Lambda}_{2,k} & \mathbf{\Lambda}_{3,k}
        \end{bmatrix}\in \text{C}_{2n+1}\, \forall k \in [K]\\
        \text{(C8)} \quad & \bm{\lambda}_{2,k}\in\mathbb{R}_n^+, \mathbf{\Lambda}_{2,k}\in\mathbb{R}_{n\times n}^+, \mathbf{\Lambda}_{3,k}\in\mathbb{S}_n\cap\mathbb{R}_{n\times n}^+ \quad \forall k \in [K]
    \end{aligned}\right\}.
\end{align}
\fi
We get constraint \eqref{C5} since:
\begin{align*}
    \sum_{k\in[K]}&\left(\lambda_{1,k}(i)+\lambda_{2,k}(i)\right)=\sum_{k\in[K]}\left(\mathbb{E}[\tilde{\xi}_i\mathbb{I}(k^*(\bm{\tilde{\xi}})=k)]+\mathbb{E}[(1-\tilde{\xi}_i)\mathbb{I}(k^*(\bm{\tilde{\xi}})=k)]\right)\\   &=\sum_{k\in[K]}\left(\mathbb{P}(\tilde{\xi}_i=1, k^* (\bm{\tilde{\xi}})=k)+\mathbb{P}(\tilde{\xi}_i=0, k^*(\bm{\tilde{\xi}})=k)\right) \text{[since } \tilde{\xi}_i \in \{0,1\} \text{]}.\\  &=\mathbb{P}(\tilde{\xi}_i=1)+\mathbb{P}(\tilde{\xi}_i=0) =1
\end{align*}
Similarly, 
\begin{align*}
    \sum_{k\in[K]} & \left(\Lambda_{1,k}(i,i)+2\Lambda_{2,k}(i,i)+\Lambda_{3,k}(i,i)\right)\\
    &=\sum_{k\in[K]}(\mathbb{E}[\tilde{\xi}_i^2\mathbb{I}(k^*(\bm{\tilde{\xi}})=k)]
    +2\mathbb{E}[\tilde{\xi}_i(1-\tilde{\xi}_i)\mathbb{I}(k^*(\bm{\tilde{\xi}})=k)]
    +\mathbb{E}[(1-\tilde{\xi}_i)^2\mathbb{I}(k^*(\bm{\tilde{\xi}})=k)])\\
    &=\sum_{k\in[K]}\left(\mathbb{P}(\tilde{\xi}_i=1,k^*(\bm{\tilde{\xi}})=k)+\mathbb{P}(\tilde{\xi}_i=0,k^*(\bm{\tilde{\xi}})=k)\right)=1.
\end{align*}
Again since $\tilde{\xi}_i$ is a Bernoulli random variable, we have $\tilde{\xi}_i^2 = \tilde{\xi}_i$ and taking expectations on both sides gives us constraint \eqref{C7}:
\begin{align*}
    \mathbb{E}[\bm{\tilde{\xi}}]=\text{diag}(\text{E}[\bm{\tilde{\xi}}\bm{\tilde{\xi}}^T]) \quad \Rightarrow \quad \sum_{k\in[K]}\bm{\lambda}_{1,k}=&\sum_{k\in[K]}\text{diag}(\mathbf{\Lambda_{1,k}}).
\end{align*}
The matrix in constraint \eqref{C8} is the second moment of the random vector $\bm{y}_k(\bm{\tilde{\xi}})$ where 
$
    \bm{y}_k(\bm{\tilde{\xi}})=\begin{pmatrix}
        1, & \bm{\tilde{\xi}},& \bm{1}-\bm{\tilde{\xi}}
    \end{pmatrix}'\mathbb{I}(k^*(\bm{\tilde{\xi}})=k)
$.
%Then,
%\begin{align*}
    %\text{E}\left[\bm{y}_k(\bm{\tilde{\xi}})\bm{y}_k(\bm{\tilde{\xi}})^T\right]=&\text{E}\left[\begin{pmatrix}
    %    1 & \bm{\tilde{\xi}}^T & (\bm{1}-\bm{\tilde{\xi}})^T\\
     %   \bm{\tilde{\xi}} & \bm{\tilde{\xi}}\bm{\tilde{\xi}}^T & \bm{\tilde{\xi}}(\bm{1}-\bm{\tilde{\xi}})^T\\
     %   (\bm{1}-\bm{\tilde{\xi}}) & (\bm{1}-\bm{\tilde{\xi}})\bm{\tilde{\xi}}^T & (\bm{1}-\bm{\tilde{\xi}})(\bm{1}-\bm{\tilde{\xi}})^T
    %\end{pmatrix}\mathbb{I}(k^*(\bm{\tilde{\xi}})=k)\right]\\
    %=&\begin{pmatrix}
     %   \lambda_{0,k} & \bm{\lambda}_{1,k}^T & \bm{\lambda}_{2,k}^T\\
     %   \bm{\lambda}_{1,k} & \mathbf{\Lambda}_{1,k} & \mathbf{\Lambda}_{2,k}^T\\
     %   \bm{\lambda}_{2,k} & \mathbf{\Lambda}_{2,k} & \mathbf{\Lambda}_{3,k}
    %\end{pmatrix}
%\end{align*}
Since $\bm{\tilde{\xi}}\in\{0,1\}^n$, $\bm{y}_k(\bm{\tilde{\xi}})$ is a nonnegative vector and $\bm{y}_k(\bm{\tilde{\xi}})\bm{y}_k(\bm{\tilde{\xi}})^T$ is a completely positive matrix. Due to the convexity of the set of all completely matrices, $\mathbb{E}\left[\bm{y}_k(\bm{\tilde{\xi}})\bm{y}_k(\bm{\tilde{\xi}})^T\right]$ is also completely positive.

\subsection{Support $\Xi = \mathbb{R}_n^+$}
\label{subsec:rnplus_support}

For the support of $\mathbb{R}_n^+$,
\begin{align}
    \mathcal{F}=\left\{(\lambda_{0,k}, \bm{\lambda}_{1,k}, \mathbf{\Lambda}_{1,k}) \text{ } \forall k\in[K]: \quad \begin{bmatrix}
        \lambda_{0,k} & \bm{\lambda}^T_{1,k} \\
        \bm{\lambda}_{1,k} & \mathbf{\Lambda}_{1,k}
    \end{bmatrix}\in \text{C}_{n+1} \tag{C9}, \ \forall k\in[K]\right\}.
\end{align}
Since $\bm{\tilde{\xi}}\in\mathbb{R}^+_n$,  the moment matrix here is a completely positive matrix due to
\begin{align*}
    \mathbb{E}\left[\begin{pmatrix}
        1 \\
        \bm{\tilde{\xi}} 
    \end{pmatrix}\begin{pmatrix}
        1 \\
        \bm{\tilde{\xi}} 
    \end{pmatrix}^T\mathbb{I}(k^*(\bm{\tilde{\xi}})=k)\right]&=\mathbb{E}\left[\begin{pmatrix}
        1 & \bm{\tilde{\xi}}^T\\
        \bm{\tilde{\xi}} & \bm{\tilde{\xi}}\bm{\tilde{\xi}}^T
    \end{pmatrix}\mathbb{I}(k^*(\bm{\tilde{\xi}})=k)\right]\\
    &=\begin{bmatrix}
        \mathbb{P}(k^*(\bm{\tilde{\xi}})=k) & \mathbb{E}\left[\bm{\tilde{\xi}}\mathbb{I}(k^*(\bm{\tilde{\xi}})=k)\right]^T \\
        \mathbb{E}\left[\bm{\tilde{\xi}}\mathbb{I}(k^*(\bm{\tilde{\xi}})=k)\right] & \mathbb{E}\left[\bm{\tilde{\xi}}\bm{\tilde{\xi}}^T\mathbb{I}(k^*(\bm{\tilde{\xi}})=k)\right]
    \end{bmatrix}\in\text{C}_{n+1}
\end{align*}

\subsection{Support $\Xi  = \mathbb{R}_n$}
\label{subsec:rn_support}

For the support of $\mathbb{R}_n$,
\begin{align}
    \mathcal{F}=\left\{(\lambda_{0,k}, \bm{\lambda}_{1,k}, \mathbf{\Lambda}_{1,k}) \text{ } \forall k\in[K]: \tag{C10} \quad \begin{bmatrix}
        \lambda_{0,k} & \bm{\lambda}^T_{1,k} \\
        \bm{\lambda}_{1,k} & \mathbf{\Lambda}_{1,k}
    \end{bmatrix}\succeq 0, \ \forall k\in[K]\right\}
\end{align}
The positive definiteness is a result of the moment matrix:
\begin{align*}
    \mathbb{E}\left[\begin{pmatrix}
        1 & \bm{\tilde{\xi}}^T\\
        \bm{\tilde{\xi}} & \bm{\tilde{\xi}}\bm{\tilde{\xi}}^T
    \end{pmatrix}\mathbb{I}(k^*(\bm{\tilde{\xi}})=k)\right]=\begin{bmatrix}
        \mathbb{P}(k^*(\bm{\tilde{\xi}})=k) & \mathbb{E}\left[\bm{\tilde{\xi}}\mathbb{I}(k^*(\bm{\tilde{\xi}})=k)\right]^T \\
    \mathbb{E}\left[\bm{\tilde{\xi}}\mathbb{I}(k^*(\bm{\tilde{\xi}})=k)\right] & \mathbb{E}\left[\bm{\tilde{\xi}}\bm{\tilde{\xi}}^T\mathbb{I}(k^*(\bm{\tilde{\xi}})=k)\right]
    \end{bmatrix}\succeq 0.
\end{align*}
For the case of $\Xi = \mathbb{R}_n$, we note that a very similar  formulation has been proposed for the objective of maximum of  linear functions of $\bm{\tilde{\xi}}$ \cite{Mishra2012}. 
%However in the existing formulation exact knowledge of first two moments is assumed whereas we assume bounds on the moments. Also we look at the more general objective of  maximum of quadratic functions.
\\
\\By construct, $value\eqref{CPP}\geq value\eqref{sup_problem}$ for all the supports where $value(\cdot)$ returns the optimal value of the program.

\section{Tightness Proofs}
\label{sec:tightnessProofs}
In this section, we demonstrate that the general conic program \eqref{CPP} is a tight formulation to \eqref{sup_problem} for all support cases. First, we consider the Bernoulli support case. We build on the results of \cite{Burer09} and adapt them to our setting to establish tightness. We subsequently extend our results to $\mathbb{R}_n^+$ and $\mathbb{R}_n$ support. 

\subsection{Tightness for  $\Xi = \{0,1\}^n$}\label{sec:bernoulli}
Consider any completely positive decomposition of the feasible solutions to \eqref{CPP}:
\begin{align*}
    \begin{pmatrix}
        \lambda_{0,k} & \bm{\lambda}_{1,k}^T & \bm{\lambda}_{2,k}^T\\
        \bm{\lambda}_{1,k} & \mathbf{\Lambda}_{1,k} & \mathbf{\Lambda}_{2,k}^T\\
        \bm{\lambda}_{2,k} & \mathbf{\Lambda}_{2,k} & \mathbf{\Lambda}_{3,k}
    \end{pmatrix} = \sum_{l=1}^{L_k}\begin{pmatrix}
        \alpha_{k,l} \\ \bm{\beta}_{k,l} \\ \bm{\eta}_{k,l}
    \end{pmatrix}\begin{pmatrix}
        \alpha_{k,l} \\ \bm{\beta}_{k,l} \\ \bm{\eta}_{k,l}
    \end{pmatrix}^T
\end{align*}
where $\alpha_{k,l}\in\mathbb{R}^+$, $\bm{\beta}_{k,l},\bm{\eta}_{k,l}\in\mathbb{R}^+_n$ $\forall l\in [L_k]$ and $\forall k \in [K]$. Without loss of generality, the above completely positive decomposition can be rewritten as
\begin{align}\label{eqn:cpp_decomposition}
    \begin{pmatrix}
        \lambda_{0,k} & \bm{\lambda}_{1,k}^T & \bm{\lambda}_{2,k}^T\\
        \bm{\lambda}_{1,k} & \mathbf{\Lambda}_{1,k} & \mathbf{\Lambda}_{2,k}^T\\
        \bm{\lambda}_{2,k} & \mathbf{\Lambda}_{2,k} & \mathbf{\Lambda}_{3,k}
    \end{pmatrix} = \sum_{t=T_{k-1}+1}^{T_k}\begin{pmatrix}
        \alpha_t\\ \bm{\beta}_t \\ \bm{\eta}_t
    \end{pmatrix}\begin{pmatrix}
        \alpha_t\\ \bm{\beta}_t \\ \bm{\eta}_t
    \end{pmatrix}^T
\end{align}
where $T_k=L_0+L_1+\dots+L_k$ for $k\in[K]$ and $L_0=0$. 
\\We tailor the properties of the CPP solutions by \cite{Burer09} to our setting.
\begin{proposition}[Burer, 2009]\label{prop:decomposition}
    For the decomposition of \eqref{eqn:cpp_decomposition}, define $\mathcal{T}^0=\{t\in[T]:\alpha_t=0\}$ and $\mathcal{T}^+=\{t\in[T]:\alpha_t>0\}$ where $T=\sum_{k\in[K]\cup\{0\}}L_k$. Then,
    \begin{enumerate}
        \item $\bm{\eta}_t/\alpha_t=\bm{1}-\bm{\beta}_t/\alpha_t$ and $\bm{\beta}_t/\alpha_t\leq\bm{1}$ for $\forall t\in \mathcal{T}^+$,
        \item $\bm{\eta}_t=\bm{\beta}_t=\bm{0}$ for $\forall t\in\mathcal{T}^0$.
    \end{enumerate}
\end{proposition}
{\it Proof.} From decomposition \eqref{eqn:cpp_decomposition}, we have for all $k \in [K]$
\begin{align*}
    &\lambda_{0,k}=\sum_{t=1 + T_k - 1}^{T_k}\alpha_t^2, \quad \bm{\lambda}_{1,k}=\sum_{t=1+T_{k-1}}^{T_k}\alpha_t\bm{\beta}_t, \quad \bm{\lambda}_{2,k}=\sum_{t=1+T_{k-2}}^{T_k}\alpha_t\bm{\eta}_t,\\
    &\mathbf{\Lambda}_{1,k}=\sum_{t=1+T_{k-1}}^{T_k}\bm{\beta}_t\bm{\beta}_t^T,\quad \mathbf{\Lambda}_{2,k}=\sum_{t=1+T_{k-1}}^{T_k}\bm{\eta}_t\bm{\beta}_t^T, \quad \text{and} \quad \mathbf{\Lambda}_{3,k}=\sum_{t=1+T_{k-1}}^{T_k}\bm{\eta}_t\bm{\eta}_t^T.
\end{align*}
Consider the constraints \eqref{C5} and \eqref{C6}:
\begin{align*}
\begin{array}{cc}
    \sum_{k\in[K]}\sum_{t=1+T_{k-1}}^{T_k}\alpha_t(\beta_t(i)+\eta_t(i))=1,\\
    \sum_{k\in[K]}\sum_{t=1+T_{k-1}}^{T_k}(\beta_t^2(i)+2\beta_t(i)\eta_t(i)+\eta_t(i)^2)=1
\end{array}
    \quad \Rightarrow \quad \begin{array}{cc}
         \sum_{t\in[T]}\alpha_t(\beta_t(i)+\eta_t(i))=1,\\
    \sum_{t\in[T]}(\beta_t(i)+\eta_t(i))^2=1
    \end{array}
\end{align*}
for $\forall i \in [n]$ where $T=\sum_{k\in[K]\cup\{0\}}L_k$.
\\From the constraint \eqref{C4},
\begin{align*}
    \left(\sum_{t\in[T]}\alpha_t(\beta_t(i)+\eta_t(i))\right)^2=\sum_{t\in[T]}\alpha_t^2\sum_{t\in[T]}(\beta_t(i)+\eta_t(i))^2
\end{align*}
By the equality conditions of Cauchy-Schwartz inequality, $\exists\delta_i$ such that $\delta_i\alpha_t=\beta_t(i)+\eta_t(i)$ for $\forall i\in[n]$ and $\forall t\in[T]$. In fact, from the results in \cite{Burer09}$, \delta_i=1$ for $\forall i \in [n]$. Then,
\begin{align*}
    \beta_t(i)+\eta_t(i)=0 \quad \forall t\in {\mathcal{T}^0} \quad &\Rightarrow \quad \beta_t(i)=\eta_t(i)=0
\end{align*}
as $\beta_t(i),\eta_t(i)\geq 0$. Hence, (ii) holds. Furthermore,
\begin{align*}
    \frac{\beta_t(i)}{\alpha_t}+\frac{\eta_t(i)}{\alpha_t}=1 \quad \forall t\in {\mathcal{T}^+} \quad &\Rightarrow \quad \frac{\eta_t(i)}{\alpha_t}=1-\frac{\beta_t(i)}{\alpha_t}\geq 0\\
    &\Rightarrow 0\leq \frac{\beta_t(i)}{\alpha_t} \leq 1.
\end{align*}
Hence, (i) holds.
\qed
Let $p_t=\alpha_t^2$ and $\bm{\zeta}_t=\bm{\beta}_t/\alpha_t$ for $\forall t\in{\mathcal{T}^+}$, then decomposition \eqref{eqn:cpp_decomposition} can be rewritten as
\begin{align}\label{eqn:cpp_decomposition_2}
    \begin{pmatrix}
        \lambda_{0,k} & \bm{\lambda}_{1,k}^T & \bm{\lambda}_{2,k}^T\\
        \bm{\lambda}_{1,k} & \mathbf{\Lambda}_{1,k} & \mathbf{\Lambda}_{2,k}^T\\
        \bm{\lambda}_{2,k} & \mathbf{\Lambda}_{2,k} & \mathbf{\Lambda}_{3,k}
    \end{pmatrix} = \sum_{t\in\mathcal{T}_k^+}p_t\begin{pmatrix}
        1\\ \bm{\zeta}_t \\ \bm{1}-\bm{\zeta}_t
    \end{pmatrix}\begin{pmatrix}
        1\\ \bm{\zeta}_t \\ \bm{1}-\bm{\zeta}_t
    \end{pmatrix}^T
\end{align}
where $\mathcal{T}_k^+=\{t\in [T_k]/[T_{k-1}]: \alpha_t>0\}$, $\forall k\in[K]$; $\sum_{t\in\mathcal{T}^+}p_t=1$ for $\mathcal{T}^+=\bigcup_{k\in[K]}\mathcal{T}_k^+$; and $\bm{0}\leq \bm{\zeta}_t\leq \bm{1}$, $\forall t\in\mathcal{T}^+_k$, $\forall k\in[K]$.
\begin{proposition}[Burer, 2009]\label{prop:binary}
    According to the decomposition \eqref{eqn:cpp_decomposition_2}, $\zeta_t(i)\in\{0,1\}$ for $\forall i\in[n]$, $\forall t\in\mathcal{T}_+$.
\end{proposition}
{\it Proof.} We have 
\begin{align*}
    \bm{\lambda}_{1,k}=\sum_{t\in\mathcal{T}^+_k}p_t\bm{\zeta}_t \quad \text{and} \quad \mathbf{\Lambda}_{1,k}=\sum_{t\in\mathcal{T}^+_k}p_t\bm{\zeta}_t\bm{\zeta}_t^T.
\end{align*}
By Proposition \ref{prop:decomposition}(ii), $\zeta_t(i)\in[0,1]$ for $\forall i\in [n]$ and $\forall t\in\mathcal{T}_k^+, \forall k\in[K]$. Then, $\zeta_t(i)^2\leq\zeta_t(i)$. 
\\From  constraint \eqref{C7}:
\begin{align*}
    \sum_{t\in\mathcal{T}^+}p_t\zeta_t(i)=\sum_{t\in\mathcal{T}^+}p_t\zeta_t(i)^2 \quad &\Rightarrow \quad \sum_{t\in\mathcal{T}^+}p_t(\zeta_t(i)-\zeta_t(i)^2)=0.
\end{align*}
Since $p_t>0$ and $\zeta_t(i)-\zeta_t(i)^2\geq 0$ for $\forall t\in\mathcal{T}^+$, we have $\zeta_t(i)-\zeta_t(i)^2=0$. Hence, $\zeta_t(i)\in\{0,1\}$.
\qed
Now, we present our main result on the equivalence of \eqref{CPP} and \eqref{sup_problem}.
\begin{theorem}\label{thm:bernoulli}
    Suppose $\mathcal{P}$ is a non-empty set, then $value\eqref{sup_problem}=value\eqref{CPP}$ for  $\Xi = \{0,1\}^n$.
\end{theorem}
{\it Proof.} Due to the interpretation for the decision variables and constraints provided in \Cref{subsec:01support}, we know that $value\eqref{CPP}\geq value\eqref{sup_problem}$. Now it suffices to show that $value\eqref{CPP}\leq value\eqref{sup_problem}$. Consider the optimal solution to \eqref{CPP}, denoted as 
$
    \left\{(\lambda_{0,k}^\ast, \bm{\lambda}_{1,k}^\ast, \bm{\lambda}_{2,k}^\ast, \mathbf{\Lambda}_{1,k}^\ast,\mathbf{\Lambda}_{2,k}^\ast, \mathbf{\Lambda}_{3,k}^\ast)\right\}_{k\in[K]}.
$
The optimal solution has the following completely positive decomposition by Propositions \ref{prop:decomposition} and \ref{prop:binary},
\begin{align*}
    \begin{bmatrix}
        \lambda_{0,k}^\ast & {\bm{\lambda}^\ast_{1,k}}^T & {\bm{\lambda}^\ast_{2,k}}^T \\
        {\bm{\lambda}^\ast_{1,k}} & {\mathbf{\Lambda}_{1,k}^\ast} & {\mathbf{\Lambda}_{2,k}^\ast}^T \\
        {\bm{\lambda}^\ast_{2,k}} & {\mathbf{\Lambda}_{2,k}^\ast} & {\mathbf{\Lambda}_{3,k}^\ast}
    \end{bmatrix}=&\sum_{t\in\mathcal{T}_k^{+}} p_t^\ast\begin{pmatrix}
        1 \\ \bm{\zeta}_t \\ \bm{1}-\bm{\zeta}_t
    \end{pmatrix}\begin{pmatrix}
        1 \\ \bm{\zeta}_t \\ \bm{1}-\bm{\zeta}_t
    \end{pmatrix}^T
\end{align*}
where $\bm{\zeta}_t\in\{0,1\}^n$, $\forall t \in\mathcal{T}_k^+$ and $\forall k\in[K]$. The decomposition allows for the  construction of a feasible distribution.
%with the desired support since,
%\begin{align*}
%    \mathbb{P}^\ast\left(\bm{\tilde \xi}=\bm{\zeta}_t, k^*(\bm{\tilde{\xi}})=k\right)&=p^*_t \quad \text{for} \quad \forall t \in\mathcal{T}_k^+, \forall k\in[K] \text{ and},\\
%    \sum_{k\in[K]} \sum_{t\in\mathcal{T}^+_k} \mathbb{P}^\ast\left(\bm{\tilde \xi}=\bm{\zeta}_t, k^*(\bm{\tilde{\xi}})=k\right)=&\sum_{t\in\mathcal{T}^+} p^\ast_t 
%    =1 \quad \text{(due to constraint (C3)).}
%\end{align*}
In particular a distribution $\mathbb{P}^*$ over $\bm{\tilde{\xi}}$ can be constructed as follows:
\begin{enumerate}
    \item Generate a random variable $\tilde{z} = k$ with probability $\lambda_{0,k}^*$.
    \item Given the realization of $\tilde{z} =k$, generate $\bm{\tilde{\xi}}$ using the following conditional probability distribution, $\mathbb{P}(\bm{\tilde{\xi}} = \bm{\zeta_t}| \tilde{z} = k) = p_t^*/\lambda_{0,k}^\ast$ $ \forall t \in \mathcal{T}_k^+$.
\end{enumerate}
$\mathbb{P}^*$ is a valid distribution since,
$
    \sum_{k \in K} \lambda_{0,k}^* \sum_{t \in \mathcal{T}_k^+} p^*_t / \lambda^*_{0,k} = \sum_{t \in \mathcal{T}^+} p^*_t = 1.
$
The constructed distribution $\mathbb{P}^* \in \mathcal{P}$ since,
\begin{align*}
    \mathbb{E}_{\mathbb{P}^*}[\bm{\tilde{\xi}}]= \sum_{k \in [K]}\lambda_{0,k}^* \sum_{t \in \mathcal{T}_k^+} p^*_t \bm{\zeta_t}/ \lambda_{0,k} = \sum_{k \in [K]} \sum_{t \in \mathcal{T}_k^+} p^*_t \bm{\zeta}_t = \sum_{k \in [K]}\bm{\lambda}^*_{1,k} 
\end{align*} 
Constraint \eqref{C1} in  \eqref{CPP} ensures that $\sum_{k \in [K]}\lambda^*_{1,k}(i) \ = \mu_i \forall i \in [n] $.
%The mean of the marginal distribution of $\bm{\tilde \xi}\mathbb{I}(k^*(\bm{\tilde{\xi}})=k)$ is
%\begin{align*}
   % \text{E}_{\mathbb{P}^\ast}\left[\bm{\tilde \xi}\mathbb{I}(k^\ast(\bm{\xi})=k)\right]=&\sum_{t\in\mathcal{T}_k^+}\bm{\zeta}_t \mathbb{P}^\ast\left(\bm{\tilde \xi}=\bm{\zeta}_t, k^*(\bm{\tilde{\xi}})=k\right)=\sum_{t\in\mathcal{T}_k^+}p^\ast_t\bm{\zeta}_t
%\end{align*}
%and satisfies $
%    \text{E}_{\mathbb{P}^\ast}[\bm{\tilde \xi}]=\sum_{t\in\mathcal{T}^+} p^\ast_t\bm{\zeta}_t\in [\bm{\underline{\gamma}}, \bm{\bar{\gamma}}]$.
Also, 
\begin{align*}
    \mathbb{E}_{\mathbb{P}^*}[\bm{\tilde{\xi}}\bm{\tilde{\xi}}^{\top}]= \sum_{k \in [K]} \lambda_{0,k}^* \sum_{t \in \mathcal{T}_k^+} p_t \bm{\zeta}_t \bm{\zeta_t}^{\top}/\lambda_{0,k} = \sum_{k \in [K]} \Lambda_{1,k}^*
\end{align*}
From  constraint \eqref{C2} in \eqref{CPP}, we have $\sum_{k \in [K]} \Lambda_{1,k}^* (i,i) = \Sigma(i,i)$ for all $i \in [n]$ and from constraint \eqref{C3} in \eqref{CPP}, we have $\sum_{k \in [K]} \Lambda_{1,k}^* (i,j) \geq \Sigma(i,j)$ for all $(i,j) \in \binom{[n]}{2}$.
%Similarly, the second moment matrix of $\bm{\tilde \xi}\mathbb{I}(k^*(\bm{\tilde{\xi}})=k)$ is
%\begin{align*}  \text{E}_{\mathbb{P}^\ast}\left[\bm{\tilde \xi}\bm{\hat\xi}^T\mathbb{I}(k^*(\bm{\tilde{\xi}})=k)\right]=&\sum_{t\in\mathcal{T}_k^+}\bm{\zeta}_t\bm{\zeta}_t^TP^\ast\left(\bm{\tilde \xi}=\bm{\zeta}_t, k^*(\bm{\tilde{\xi}})=k\right)=\sum_{t\in\mathcal{T}_k^+}p_t^\ast\bm{\zeta}_t\bm{\zeta}_t^T
%\end{align*}
%and satisfies
%\begin{align*}
%    \text{E}_{\mathbb{P}^\ast}\left[\tilde{\xi}_i{\hat\xi}_j\mathbb{I}(k^*(\bm{\tilde{\xi}})=k)\right]=\sum_{t\in\mathcal{T}^+}p_t\bm{\zeta}_t\bm{\zeta}_t^T\geq\Sigma_{ij} \quad \forall (i,j)\in[n]_2.
%\end{align*}
Thus, the constructed distribution is a feasible solution to \eqref{sup_problem}, and further,
\begin{align*}
    value\eqref{sup_problem}& \geq\text{E}_{\mathbb{P}^*}\left[\max_{k\in[K]}\left\{\bm{\tilde \xi}^T\bm{A}_k\bm{\tilde \xi}+\bm{\tilde \xi}^T\bm{b}_k+c_k\right\}\right]\\
    & \geq \sum_{k\in[K]} \mathbb{P}^*(\tilde{z} = k)\left(\mathbb{E}_{\mathbb{P}^*}\left[\bm{\tilde \xi}\bm{\tilde \xi}^T| \tilde{z} = k \right]\bullet\mathbf{A}_k+\mathbb{E}_{\mathbb{P}^*}\left[\bm{\tilde \xi}| \tilde{z} = k\right]^T\bm{b}_k
    +c_k\right)\\
=&\sum_{k\in[K]}\lambda_{0,k}\left[\left(\sum_{t\in\mathcal{T}_k^+}\frac{p_t^\ast\bm{\zeta}_t\bm{\zeta}_t^T}{\lambda_{0,k}^\ast}\right)\bullet\mathbf{A}_k+\left(\sum_{t\in\mathcal{T}_k^+}\frac{p_t^\ast\bm{\zeta}_t}{\lambda_{0,k}^\ast}\right)^T\bm{b}_k+c_k\right]\\
&=\sum_{k\in[K]}\left[\bm{\Lambda}_{1,k}^*\bullet\mathbf{A}_k+ (\bm{\lambda}_{1,k}^*)^T\bm{b}_k+ \lambda_{0,k}^\ast c_k\right]\\
    &= value\eqref{CPP}.
\end{align*}
The first inequality is true because $\mathbb{P}^*$ is a feasible distribution, the second inequality arises since we use a specific value of $k$ instead of the maximizer.
%\begin{align*}
%    value\eqref{sup_problem}\geq& \text{E}_{\mathbb{P}^*}\left[\max_{k\in[K]}\left\{\bm{\tilde \xi}^T\bm{A}_k\bm{\tilde \xi}+\bm{\tilde \xi}^T\bm{b}_k+c_k\right\}\right]\\
    %=&\sum_{k\in[K]}\Bigg[\text{E}_{\mathbb{P}^*}\left[\bm{\tilde \xi}\bm{\tilde \xi}^T\mathbb{I}(k^\ast(\bm{\xi})=k)\right]\bullet\mathbf{A}_k+\text{E}_{\mathbb{P}^*}\left[\bm{\tilde \xi}\mathbb{I}(k^\ast(\bm{\xi})=k)\right]^T\bm{b}_k\\
%    &+P^\ast(k^\ast(\bm{\xi})=k)c_k\Bigg]\\
    %=&\sum_{k\in[K]}\left[\left(\sum_{t\in\mathcal{T}_k^+}p_t^\ast\bm{\zeta}_t\bm{\zeta}_t^T\right)\bullet\mathbf{A}_k+\left(\sum_{t\in\mathcal{T}_k^+}p_t^\ast\bm{\zeta}_t\right)^T\bm{b}_k+\left(\sum_{t\in\mathcal{T}_k^+}p_t^\ast\right)c_k\right]\\
 %   =&value\eqref{CPP}.
%\end{align*}
Therefore since $value\eqref{sup_problem}\leq value\eqref{CPP}$ by construction, together with this proof $value\eqref{sup_problem}=value\eqref{CPP}$.
\qed

\subsection{Tightness for supports $\mathbb{R}^+_n$ and $\mathbb{R}_n$}
Similar to the case of $\Xi = \{0,1\}^n$, consider any decomposition of the feasible solutions to \eqref{CPP}:
\begin{align*}
    \begin{pmatrix}
        \lambda_{0,k} & \bm{\lambda}_{1,k}^T & \\
        \bm{\lambda}_{1,k} & \mathbf{\Lambda}_{1,k} 
    \end{pmatrix} = \sum_{l=1}^{L_k}\begin{pmatrix}
        \alpha_{k,l} \\ \bm{\beta}_{k,l} 
    \end{pmatrix}\begin{pmatrix}
        \alpha_{k,l} \\ \bm{\beta}_{k,l} 
    \end{pmatrix}^T
\end{align*}
where $\forall k \in [K], \forall l\in [L_k], \alpha_{k,l}\in\mathbb{R}^+$, $\bm{\beta}_{k,l} \in \mathbb{R}_n^+ $ for $\Xi = \mathbb{R}_n^+$ and  and $\bm{\beta}_{k,l} \in \mathbb{R}_n $ for $\Xi = \mathbb{R}_n$. Without loss of generality, the above decomposition can be rewritten as
\begin{align}\label{eqn:gen_decomp}
    \begin{bmatrix}
        \lambda_{0,k} & {\bm{\lambda}_{1,k}}^T\\
        {\bm{\lambda}_{1,k}} & {\mathbf{\Lambda}_{1,k}}
    \end{bmatrix}=\sum_{t\in\mathcal{T}_k}\begin{pmatrix}
        \alpha_t \\
        \bm{\beta}_t
    \end{pmatrix}\begin{pmatrix}
        \alpha_t \\
        \bm{\beta}_t
    \end{pmatrix}^T
\end{align}
where $(\alpha_t, \bm{\beta}_t)\in\mathbb{R}^+\times\mathbb{R}^+_n$ for the support $\mathbb{R}_n^+$, $(\alpha_t, \bm{\beta}_t)\in\mathbb{R}^+\times\mathbb{R}_n^+$ for the support $\mathbb{R}_n$. Let $T_k=L_0+L_1+\dots+L_k$ for $k\in[K]$ and $L_0=0$, then $\mathcal{T}_k = \{ T_{k-1} +1, \ldots, T_{k} \}$.  %$\mathcal{T}_k$ is defined as in Section \ref{sec:bernoulli}. 
Furthermore, we can rewrite the decomposition \eqref{eqn:gen_decomp} as
\begin{align}\label{eqn:gen_sep_decomp}
    \begin{bmatrix}
        \lambda_{0,k} & {\bm{\lambda}_{1,k}}^T\\
        {\bm{\lambda}_{1,k}} & {\mathbf{\Lambda}_{1,k}}
    \end{bmatrix}=\sum_{t\in\mathcal{T}_k^+}{(\alpha_t)}^2\begin{pmatrix}
        1 \\
        \bm{\beta}_t/\alpha_t
    \end{pmatrix}\begin{pmatrix}
        1 \\
        \bm{\beta}_t/\alpha_t
    \end{pmatrix}^T + \sum_{t\in\mathcal{T}_k^0}\begin{pmatrix}
        0 \\
        \bm{\beta}_t
    \end{pmatrix}\begin{pmatrix}
        0 \\
        \bm{\beta}_t
    \end{pmatrix}^T
\end{align}
where $\mathcal{T}^+_k$ and $\mathcal{T}_k^0$ are defined (similar to \Cref{prop:decomposition}) as follows, $\mathcal{T}^0=\{t\in \mathcal{T}_k:\alpha_t=0\}$ and $\mathcal{T}^+=\{t\in \mathcal{T}_k:\alpha_t>0\}$. 
We will make use of this decomposition in our proof of tightness which follows.
\begin{theorem}\label{thm:R_n_+_and_R_n}
    Suppose $\mathcal{P}$ is a non-empty set, then $value\eqref{sup_problem}=value\eqref{CPP}$ for supports $\mathbb{R}_n^+$ and $\mathbb{R}_n$.
\end{theorem}
{\it Proof. }Similar to the proof of Theorem \ref{thm:bernoulli}, it is sufficient to show that $value\eqref{CPP}\leq value\eqref{sup_problem}$ for both support cases. The proof of Theorem \ref{thm:R_n_+_and_R_n} relies on constructing a sequence of random vectors and feasible solutions with respect to the decomposition \eqref{eqn:gen_sep_decomp}.
 Consider the optimal solution to \eqref{CPP} denoted as $\left\{\left(\lambda_{0,k}^\ast, \bm{\lambda}_{1,k}^\ast, \mathbf{\Lambda}^\ast_{1,k}\right)\right\}_{k\in[K]}$. The optimal solution can be written in the decomposition \eqref{eqn:gen_sep_decomp} form as follows:
\begin{align*}
    \begin{bmatrix}
        \lambda_{0,k}^\ast & {\bm{\lambda}_{1,k}^\ast}^T\\
        {\bm{\lambda}_{1,k}^\ast} & {\mathbf{\Lambda}_{1,k}^\ast}
    \end{bmatrix}=\sum_{t\in\mathcal{T}_k^+}{(\alpha_t^\ast)}^2\begin{pmatrix}
        1 \\
        \bm{\beta}_t^\ast/\alpha_t^\ast
    \end{pmatrix}\begin{pmatrix}
        1 \\
        \bm{\beta}_t^\ast/\alpha_t^\ast
    \end{pmatrix}^T + \sum_{t\in\mathcal{T}_k^0}\begin{pmatrix}
        0 \\
        \bm{\beta}_t^\ast
    \end{pmatrix}\begin{pmatrix}
        0 \\
        \bm{\beta}_t^\ast
    \end{pmatrix}^T
\end{align*}
{\it Case (1):  $\Xi = \mathbb{R}^+_n$}. We construct a sequence of random vectors so that the limiting distribution of the sequence is feasible to \eqref{sup_problem}. Let $\epsilon\in(0,1)$, then construct the distribution as follows.
\begin{enumerate}
\item Generate a random variable $\tilde{z} = k$ with probability $\lambda_{0,k}^\ast$.
\item Given the realization of $\tilde{z} = k$, generate $\bm{\tilde{\xi}}_\epsilon$ using the following conditional distribution.
\begin{align*} \mathbb{P}^\ast\left(\tilde{\bm{\xi}}_\epsilon=\bm{\beta}_t^\ast/\alpha^\ast_t \middle| \tilde{z}=k\right)&=(1-\epsilon^2)(\alpha_t^*)^2/\lambda_{0,k}^* \quad \forall t\in\mathcal{T}_k^+ \text{ and } \forall k\in[K],\\
\mathbb{P}^\ast\left(\tilde{\bm{\xi}}_\epsilon=\frac{\bm{\beta}_t\sqrt{|\mathcal{T}^0|}}{\epsilon} \middle|\tilde{z}=k\right)&=\frac{\epsilon^2}{|\mathcal{T}^0|\lambda^*_{0,k}} \quad \forall t\in\mathcal{T}_k^0 \text{ and } \forall k\in[K]
\end{align*}
where $\mathcal{T}^0=\bigcup_{k\in[K]}\mathcal{T}_k^0$. 
\end{enumerate}
This is a valid distribution as
\begin{align*}
    \sum_{k\in[K]}\lambda_{0,k}^*\left[\sum_{t\in\mathcal{T}_k^+}\mathbb{P}^\ast\Bigg(\tilde{\bm{\xi}}_\epsilon=\frac{\bm{\beta}_t^\ast}{\alpha^\ast_t} \middle| \tilde{z}=k \right)&+\sum_{t\in\mathcal{T}_k^0}\mathbb{P}^\ast\left(\tilde{\bm{\xi}}_\epsilon=\frac{\bm{\beta}_t\sqrt{|\mathcal{T}^0|}}{\epsilon} \middle| \tilde{z}=k\right)\Bigg]\\
    =&(1-\epsilon^2)\sum_{t\in\mathcal{T}^+}(\alpha_t^\ast)^2+\epsilon^2\sum_{t\in\mathcal{T}^0}\frac{1}{|\mathcal{T}^0|}\\
    =&(1-\epsilon^2)+\epsilon^2\frac{|\mathcal{T}^0|}{|\mathcal{T}^0|} \quad \text{(due to \eqref{C3})}\\
    =&1.
\end{align*}
The mean of the conditional distribution of $\bm{\tilde \xi}_\epsilon$ given $\tilde{z} = k$ is
\begin{align*}
    \lim_{\epsilon\downarrow0} \mathbb{E}_{\mathbb{P}^\ast}\left[\bm{\tilde \xi}_\epsilon \middle|\tilde{z} = k \right]=&\lim_{\epsilon\downarrow0}\sum_{t\in\mathcal{T}_k^+}\frac{\bm{\beta}_t^\ast}{\alpha_t^\ast}P^\ast\Bigg(\hat{\bm{\xi}}_\epsilon=\frac{\bm{\beta}_t^\ast}{\alpha^\ast_t} | \tilde{z} = k\Bigg)\\
    &+\lim_{\epsilon\downarrow0}\sum_{t\in\mathcal{T}_k^0}\frac{\bm{\beta}_t^\ast\sqrt{|\mathcal{T}^0|}}{\epsilon}P^\ast\left(\hat{\bm{\xi}}_\epsilon=\frac{\bm{\beta}_t\sqrt{|\mathcal{T}^0|}}{\epsilon}\middle| \tilde{z} = k\right)\\
    =&\lim_{\epsilon\downarrow0}\sum_{t\in\mathcal{T}_k^+}\frac{\bm{\beta}_t^\ast}{\alpha_t^\ast \lambda_{0,k}^*}(1-\epsilon^2)(\alpha_t^\ast)^2+\lim_{\epsilon\downarrow0}\sum_{t\in\mathcal{T}_k^0}\frac{\bm{\beta}_t^\ast\sqrt{|\mathcal{T}^0|}}{\epsilon}\frac{\epsilon^2}{|\mathcal{T}^0|\lambda_{0,k}^*}\\
    =&\sum_{t\in\mathcal{T}_k^+}\alpha_t^\ast\bm{\beta}_t^\ast/\lambda_{0,k}^*
    \end{align*}
    and satisfies
    \begin{align*}
        \lim_{\epsilon\downarrow0}\text{E}_{\mathbb{P}^\ast}[\bm{\tilde \xi}_\epsilon]=\sum_{k\in[K]}\lambda_{0,k}^*\sum_{t\in\mathcal{T}_k^+}\alpha_t^\ast\bm{\beta}_t^\ast/\lambda_{0,k}^* \quad \Rightarrow \quad \sum_{t\in\mathcal{T}^+}\alpha_t^\ast\bm{\beta}_t^\ast = \bm{\mu}.
%        \in[\bm{\underline{\gamma}}, \bm{\bar \gamma}].
    \end{align*}
    Similarly, the conditional second moment matrix is
    \begin{align*}
        \lim_{\epsilon\downarrow0} \text{E}_{\mathbb{P}^\ast}\left[\bm{\tilde \xi}_\epsilon\bm{\tilde \xi}_\epsilon^T| \tilde{z} = k\right]=&\lim_{\epsilon\downarrow0}\sum_{t\in\mathcal{T}_k^+}\frac{\bm{\beta}_t^\ast}{\alpha_t^\ast}\frac{\bm{\beta}_t^\ast}{\alpha_t^\ast}^T \frac{(1-\epsilon^2)(\alpha^\ast_t)^2}{\lambda_{0,k}^*}\\
    &+\lim_{\epsilon\downarrow0}\sum_{t\in\mathcal{T}_k^0}\frac{\bm{\beta}_t^\ast\sqrt{|\mathcal{T}^0|}}{\epsilon}\frac{{\bm{\beta}_t^\ast}^T\sqrt{|\mathcal{T}^0|}}{\epsilon}\frac{\epsilon^2}{|\mathcal{T}^0|\lambda_{0,k}^*}\\
    =&\sum_{t\in\mathcal{T}_k}{\bm{\beta}_t^\ast}{\bm{\beta}_t^\ast}^T/\lambda_{0,k}^*
    \end{align*}
    and satisfies
    \begin{align*}
        \lim_{\epsilon\downarrow0} \text{E}_{\mathbb{P}^\ast}\left[{\tilde \xi}_\epsilon(i){\tilde \xi}_\epsilon(i)\right]=\sum_{k\in[K]}\lambda_{0,k}^*\sum_{t\in\mathcal{T}_k}{{\beta}_t^\ast(i)}{{\beta}_t^\ast}(i)/\lambda_{0,k}^* \quad \Rightarrow \quad \sum_{t\in\mathcal{T}}{{\beta}_t^\ast(i)}{{\beta}_t^\ast}(i)=\Sigma_{ii} 
    \end{align*}
    for $\forall i\in[n]$ and
    \begin{align*}
        \lim_{\epsilon\downarrow0} \text{E}_{\mathbb{P}^\ast}\left[{\tilde \xi}_\epsilon(i){\tilde \xi}_\epsilon(j)\right]=\sum_{k\in[K]}\lambda_{0,k}^*\sum_{t\in\mathcal{T}_k}{{\beta}_t^\ast(i)}{{\beta}_t^\ast}(j)/\lambda_{0,k}^* \quad \Rightarrow \quad \sum_{t\in\mathcal{T}}{{\beta}_t^\ast(i)}{{\beta}_t^\ast}(j)\geq\Sigma_{ij}
    \end{align*}
    for $\forall(i,j) \in \binom{[n]}{2}$. Hence, the constructed sequence of random vectors are a feasible solution to the original problem \eqref{sup_problem}. Furthermore,
    \begin{align*}
        value\eqref{sup_problem}\geq& \lim_{\epsilon\downarrow0} \text{E}_{\mathbb{P}^\ast}\left[\max_{k\in[K]}\left\{\bm{\tilde{\xi}}_\epsilon^T\mathbf{A}_k\bm{\tilde{\xi}}_\epsilon+\bm{\tilde{\xi}}_\epsilon^T\bm{b}_k+c_k\right\}\right]\\
        \geq &\lim_{\epsilon\downarrow0} \sum_{k\in[K]}\mathbb{P}^*(\tilde{z} = k)\text{E}_{\mathbb{P}^\ast}\left[\bm{\tilde{\xi}}_\epsilon\bm{\tilde{\xi}}_\epsilon^T| \tilde{z} =k\right]\bullet\mathbf{A}_k\\
        &+\lim_{\epsilon\downarrow0}\sum_{k\in[K]}\mathbb{P}^*(\tilde{z} = k)\text{E}_{\mathbb{P}^\ast}\left[\bm{\tilde{\xi}}_\epsilon|\tilde{z} =k\right]^T\bm{b}_k+\lim_{\epsilon\downarrow0}\sum_{k\in[K]}\mathbb{P}^*(\tilde{z}=k)c_k\\
        =&\sum_{k\in[K]}\lambda_{0,k}^*\Bigg[\left(\sum_{t\in\mathcal{T}_k}\frac{{\bm{\beta}^\ast_t}{\bm{\beta}^\ast_t}^T}{\lambda_{0,k}^*}\right)\bullet\mathbf{A}_k+\left(\sum_{t\in\mathcal{T}^+_k} \frac{\alpha^\ast_t\bm{\beta}_t^\ast}{\lambda_{0,k}^*}\right)^T\bm{b}_k+c_k\Bigg]\\
        =&value\eqref{CPP}.
    \end{align*}
    The first inequality is true because $\mathbb{P}^*$ is a feasible distribution in the limit, the second inequality arises since we use a specific value of $k$ instead of the maximizer.\\
{\it Case (2):  $\Xi = \mathbb{R}_n$}. Here, we construct the following distribution $\mathbb{P}^*$ over $\bm{\tilde{\xi}}$.
    \begin{enumerate}
    \item Generate a random variable $\tilde{z} = k$ with probability $\lambda_{0,k}^*$.
    \item Given the realization of $\tilde{z} =k$, generate $\bm{\tilde{\xi}}$ from a normal distribution with mean $\bm{\lambda}^\ast_{1,k}/\lambda_{0,k}^\ast$ and second moment  $\mathbf{\Lambda}_{1,k}^\ast/\lambda_{0,k}^\ast$ 
\end{enumerate}
The distribution is a mixture of normal distributions and hence valid. Since,
\begin{align*}
        \mathbb{E}_{\mathbb{P}^\ast}[\bm{\tilde \xi}]=&\sum_{k\in[K]} \mathbb{E}_{\mathbb{P}^\ast}\left[\bm{\tilde \xi}|\tilde{z}=k\right]\mathbb{P}^\ast(\tilde{z}=k)
        =\sum_{k\in[K]} \frac{\bm{\lambda}^\ast_{1,k}}{\lambda_{0,k}^\ast}\lambda_{0,k}^\ast \quad \Rightarrow \sum_{k\in[K]}\bm{\lambda}^\ast_{1,k}  = \bm{\mu}
    \end{align*}
     \begin{align*}
        \mathbb{E}_{\mathbb{P}^\ast}[\tilde{\xi}_i\tilde{\xi}_j]&=\sum_{k\in[K]} \mathbb{E}_{\mathbb{P}^\ast}\left[{\tilde \xi}_i{\tilde \xi}_j^T|\tilde{z}=k\right]\mathbb{P}^\ast(\tilde{z}=k)=\sum_{k\in[K]}\frac{{\Lambda}_{1,k}^\ast(i,j)}{\lambda_{0,k}^\ast}\lambda_{0,k}^\ast \\
        & \quad \Rightarrow \quad \sum_{k\in[K]}{\Lambda}_{1,k}^\ast(i,i)= \Sigma_{ii}      \quad \forall i\in[n], \text{ and}\\
        & \quad \quad \quad \sum_{k\in[K]}{\Lambda}_{1,k}^\ast(i,j)\geq \Sigma_{ij}      \quad \forall(i,j)\in\binom{[n]}{2}.
    \end{align*}
 %   \begin{align*}
  %      P^\ast(k^*(\bm{\tilde{\xi}})=k)=\lambda_{0,k}^\ast
  %  \end{align*}
  %  Define the conditional distribution of $\bm{\tilde \xi}$ given $k^\ast(\bm{\tilde \xi})=k$ as a normal distribution with mean
  %  \begin{align*}
  %      \text{E}_{\mathbb{P}^\ast}\left[\bm{\tilde \xi}|k^*(\bm{\tilde{\xi}})=k\right]=\frac{\bm{\lambda}^\ast_{1,k}}{\lambda_{0,k}^\ast} 
   % \end{align*}
   % and second moment matrix 
   % \begin{align*}
    %    \text{E}_{\mathbb{P}^\ast}\left[\bm{\tilde \xi}\bm{\tilde \xi}^T|k^*(\bm{\tilde{\xi}})=k\right]=\frac{\mathbf{\Lambda}_{1,k}^\ast}{\lambda_{0,k}^\ast}.
    %\end{align*}
  %  This is a valid conditional probability distribution as $\sum_{k\in[K]}P^\ast(k^*(\bm{\tilde{\xi}})=k)=1$ 
   % and satisfies
    %\begin{align*}
     %   \text{E}_{\mathbb{P}^\ast}[\bm{\tilde \xi}]=&\sum_{k\in[K]} \text{E}_{\mathbb{P}^\ast}\left[\bm{\tilde \xi}|k^*(\bm{\tilde{\xi}})=k\right]P^\ast(k^*(\bm{\tilde \xi})=k)\\
     %   =&\sum_{k\in[K]} \frac{\bm{\lambda}^\ast_{1,k}}{\lambda_{0,k}^\ast}\lambda_{0,k}^\ast \quad \Rightarrow \sum_{k\in[K]}\bm{\lambda}^\ast_{1,k} \in [\bm{\underline{\gamma}}, \bm{\bar \gamma}]
    %\end{align*}
    %and
    %\begin{align*}
     %   \text{E}_{\mathbb{P}^\ast}[\tilde{\xi}_i\tilde{\xi}_j]=&\sum_{k\in[K]} \text{E}_{\mathbb{P}^\ast}\left[{\tilde \xi}_i{\tilde \xi}_j^T|k^*(\bm{\tilde{\xi}})=k\right]P^\ast(k^*(\bm{\tilde \xi})=k) \\
        %&\sum_{k\in[K]}\frac{{\Lambda}_{1,k}^\ast(i,j)}{\lambda_{0,k}^\ast}\lambda_{0,k}^\ast \quad \Rightarrow \quad \sum_{k\in[K]}{\Lambda}_{1,k}^\ast(i,j)\geq \Sigma_{ij}
    %\end{align*}
    %for $\forall(i,j)\in[n]_2$. 
    We have $\mathbb{P}^* \in \mathcal{P}$. We also have,
    \begin{align*}
        value\eqref{sup_problem}& \geq \mathbb{E}_{\mathbb{P}^\ast}\left[\max_{k\in[K]}\left\{\bm{\tilde \xi}^T\mathbf{A}_k\bm{\tilde \xi}+\bm{\tilde \xi}^T\bm{b}_k+c_k\right\}\right]\\
        &\geq \sum_{k\in[K]}\lambda_{0,k}^\ast \left[\mathbb{E}_{\mathbb{P}^\ast}\left[\bm{\tilde{\xi}}\bm{\tilde{\xi}}^T|\tilde{z}=k\right]\bullet\mathbf{A}_k+\sum_{k\in[K]}\mathbb{E}_{\mathbb{P}^\ast}\left[\bm{\tilde{\xi}}|\tilde{z}=k\right]^T\bm{b}_k+ c_k\right]\\
        &=\sum_{k\in [K]}\left[\mathbf{\Lambda}_{1,k}^\ast\bullet\mathbf{A_k}+{\bm{\lambda}_{1,k}^\ast}^T\bm{b}_k+\lambda_{0,k}^\ast\right]=value\eqref{CPP}.
    \end{align*}
    Since $value\eqref{sup_problem}\leq value\eqref{CPP}$ by construct, $value\eqref{sup_problem}=value\eqref{CPP}$ for both supports.
\qed 

\begin{note}
In addition, all formulations extend in a natural manner to the case where the means are not specified but lie in an interval given by $\text{E}_\mathbb{P}[\tilde{\bm{\xi}}] \in [\underline{\bm{\gamma}},\overline{\bm{\gamma}}]$.
%\text{E}_\mathbb{P}[\tilde{\xi}_i^2] = \Sigma_{ii}, \forall i \in [n],
\end{note}

%{\bf Symbols} Please use standard mathematical notation for  symbols, that is, Italic for single letters that denote mathematical constants, variables, and unknown quantities and Roman/upright for numerals, operators, and punctuation, and commonly defined functions or abbreviations, e.g., cos, det, e or exp, lim, log, max, min, sin, tan, d (for derivative). In the definition of a set, use IR for the set of reals, conv for convex hull, cl for closure, {A: B} for the set of A such that B holds, and double parallel arrows for set-valued maps.  Other symbols include  := for equality by definition and   an empty square at the right margin to denote the end of a proof. Please limit the use of acronyms.

\section{Numerical Examples}
\label{sec:numericals}
The general conic problem \eqref{CPP} corresponds to a completely positive problem when the support is $\{0,1\}^n$ or $\mathbb{R}_n^+$. The main computational difficulty arises from the completely positive cone constraints which are not currently implementable on a solver. To address this, we replace these cones with the doubly nonnegative cone, yielding an SDP relaxation of \eqref{CPP}.

For the case where support is $\{0,1\}^n$, the problem \eqref{sup_problem} admits to an exact linear problem formulation with exponential growth with $n$. This formulation allows us to compute the optimal value without explicitly optimizing over the completely positive cone $\text{C}_{2n+1}$, and thus provides a natural benchmark for evaluating the performance of the SDP relaxation. Motivated by this, we focus on the numerical experiments for the $\{0,1\}^n$ case. Furthermore, we consider an ambiguity set in which the mean vector $\text{E}_\mathbb{P}[\bm{\tilde \xi}]$ lies within a given interval of $[\underline{\bm \gamma}, \bm{\bar \gamma}]$. This choice simplifies the generation of feasible instances compared to the ambiguity set \eqref{ambiguity_set}.

The purpose of this section is to assess the trade-off between the performance of the SDP relaxation and its computational time in the general $\{\mathbf{A}_k\}_{k=1}^K$ setting, across different problem sizes $n$.

\subsection{Experimental Setup}
We generate random instances of \eqref{sup_problem} with problem parameters $\{\mathbf{A}_k, \bm{b}_k, c_k\}_{k=1}^K$ and ambiguity set parameters $(\underline{\bm{\gamma}}, \bar{\bm{\gamma}}, \mathbf{\Sigma})$. The vectors $\{\bm{b}_k, c_k\}_{k=1}^K$ and $\bar{\bm{\gamma}}$ are sampled from random uniform distributions. In contrast, the generation of $\{\mathbf{A}_k\}_{k=1}^K$ and $\mathbf{\Sigma}$ requires additional care as the $\mathbf{A}_k$'s do not satisfy any specific structural properties beyond symmetricity, and $\mathbf{\Sigma}$  must be generated to ensure that problem \eqref{sup_problem} is feasible.

Since the matrices $\{\mathbf{A}_k\}_{k=1}^K$ are required to be symmetric, each matrix is generated as $\mathbf{A}_k=(\mathbf{M}+ \mathbf{M}^T)/2$, where the entries of $\mathbf{M}$ are drawn from a random uniform distribution. To present the general case, we exclude special structured $\{\mathbf{A}_k\}_{k=1}^K$ matrices. For example, matrices that contain all off-diagonal entries non-positive correspond to the quadratic pieces to  be submodular functions, while nonnegative off diagonal entries correspond to supermodular functions. Similary,  positive semidefinite matrices correspond to convex functions. All such cases are discarded.

To ensure the feasibility of problem \eqref{sup_problem}, the construction of the matrix $\mathbf{\Sigma}$ also depends on the values of $\bm{\bar \gamma}$. The diagonal entries of $\mathbf{\Sigma}$, which also coincides with $\underline{\gamma}_i$, must satisfy $\Sigma_{ii}\leq\bar{\gamma}_i$ for all $i\in[n]$. Assume that $\bar{\gamma}_i$ is sampled from $\mathcal{U}(l_i,1)$. Accordingly, we sample the diagonal entries as $\Sigma_{ii}\sim \mathcal{U}(l_i, \bar{\gamma}_i)$.  For the off-diagonal entries, we employ a conservative construction that guarantees feasibility. Since
\begin{align*}
    \Sigma_{ij}\leq \Lambda_{1,k}(i,j)=\mathbb{P}(\tilde{\xi}_i=1, \tilde{\xi}_j=1)\leq \min \left\{\mathbb{P}(\tilde{\xi}_i=1), \mathbb{P}(\tilde{\xi}_j=1)  \right\} \quad \forall j<i,
\end{align*}
to enforce strict but feasible lower second moment inequalities, we can set $\Sigma_{ij}\geq\min_{i}\{l_i\}$. The off diagonal entries must satisfy $\Sigma_{ij}\leq \min\{\bar{\gamma}_i, \bar{\gamma}_j\}$, which provides a valid but potentially loose upper bound. To ensure feasibility in a robust manner, we adopt a more conservative construction by restricting the sampling range to a smaller interval. Specifically, for all $i>j$, we generate
$\Sigma_{ij}\sim \mathcal{U}(d,f)$ where $d=\min_i\{l_i\}$ and $f=\min_i\{\bar{\gamma}_i-l_i\}+d$ for $\forall j<i, \forall i\in[n]$. We set $\Sigma_{ji} = \Sigma_{ij}$ so that $\Sigma$ is symmetric. %Symmetry is enforced by setting $\Sigma_{ji}=\Sigma_{ij}$. 
We summarize the bound generation for the uniform distributions of all parameters in Table \ref{tab:summary_ins_gen}.
\begin{table}[ht]\label{tab:summary_ins_gen}
% The table caption is above the table.
\caption{Bound generation scheme for random instance parameters. All variables are ultimately sampled from uniform distributions with bounds constructed as described.}
\vspace{.1 in}
\begin{center}
\begin{tabular}{|c|c|c|c|}
\hline
\textbf{Parameter} & 
\textbf{Lower Bound Generation} & \textbf{Upper Bound Generation}  \\
\hline
Entries of $\mathbf{M}$ & $l_M\sim \mathcal{U}(-6, -3)$ & $u_M\sim \mathcal{U}(1, 4)$ \\
Entries of $\bm{b}_k$ & $l_b\sim \mathcal{U}(-2, 2)$ & $u_b\sim \mathcal{U}(l_b+2, l_b+6)$\\
Scalars $c_k$ & $l_c\sim\mathcal{U}(-3, 3)$ & $u_c\sim\mathcal{U}(l_c+2, l_c+8)$ \\
Entries of $\bm{\bar \gamma}$ & $l_i\sim\mathcal{U}(0,0.5)$ & 1 \\
Diagonal entries of $\mathbf{\Sigma}$ ($\underline{\bm \gamma}$) & $l_i$ & $\bar{\gamma}_i$ \\
Off-diagonal entries of $\mathbf{\Sigma}$ & $\min_i\{l_i\}$ & $\min_i\{\bar{\gamma}_i-l_i\}+\min_i\{l_i\}$\\
\hline
\end{tabular}
\end{center}
\end{table}

All formulations are implemented in MATLAB using YALMIP \cite{Lofberg2004} as modeling interface and MOSEK \cite{mosek} as the SDP solver. For each problem instance, we generate random samples according to the procedure described above.

%\subsection{The Benchmark and Performance Metric}
We utilize the exponentially-sized linear program \eqref{EXP-LP} of problem \eqref{sup_problem} as the benchmark. For completeness, the full formulation of the linear program \eqref{EXP-LP} is provided in Appendix. Accordingly, the performance of the SDP relaxation is evaluated by the optimality gap, defined as
\begin{align*} 
\text{Optimality Gap}(\%) = 100\times\frac{|value(\text{SDP})-value(\text{EXP-LP})|}{|value(\text{EXP-LP})|}
\end{align*}

\subsection{Numerical Results} 
We conduct two sets of numerical experiments to evaluate the performance of the proposed SDP relaxation. In the first experiment, we consider a problem size of $n\in\{5, 6, 7, 8, 9, 10\}$ with $K=5$. The range of $n$ is limited due to the exponential growth of the \eqref{EXP-LP}. For each value of $n$, we generate 100 independent instances and solve the benchmark \eqref{EXP-LP} and the SDP relaxation of \eqref{CPP} for $\Xi = \{0,1\}^n$. The results on the optimality gap and computational time for both models are provided in Table \ref{tab:performance}. The trade-off between the performance of the SDP relaxation and the computational efficiency is clearly illustrated. In particular, the SDP relaxation maintains a relatively mild increase in the average optimality gap as $n$ increases, while offering significant computational advantages compared to \eqref{EXP-LP}. This shows that the SDP relaxations for  completely positive formulations  provide a practical approach to solve instances with general $\{\mathbf{A}_k\}_{k=1}^K$ structures.

\begin{table}[h]\label{tab:performance}
% The table caption is above the table.
\caption{Summary of the optimality gap of the SDP relaxation and computational time comparison with the EXP-LP benchmark}
\vspace{.1 in}
\begin{center}
\begin{tabular}{|c|ccc|ccc|ccc|}
\hline
& \multicolumn{3}{c|}{Opt Gap (\%)} 
& \multicolumn{3}{c|}{Time (SDP)} 
& \multicolumn{3}{c|}{Time (EXP-LP)} \\
$n$ 
& Avg & Min & Max 
& Avg & Min & Max 
& Avg & Min & Max \\
\hline
5 & 1.27 & 0.01 & 16.82 & 0.27 & 0.13 & 0.83 & 0.21 & 0.14 & 0.65 \\
6 & 1.43 & 0.05 & 20.02 & 0.34 & 0.21 & 0.96 & 0.49 & 0.35 & 0.60\\
7 & 5.09 & 0.02 & 131.54 & 0.43 & 0.28 & 1.07 & 1.32 & 0.78 & 1.71\\
8 & 3.38 & 0.03 & 77.69 & 0.48 & 0.31 & 0.89 & 3.68 & 2.23  & 8.60\\
9 & 4.87 & 0.06 & 67.84 & 0.63  & 0.42  & 1.38 & 10.22 & 9.45 & 11.05\\
10 & 8.69 & 0.30 & 418.87 & 0.80 & 0.34 & 3.64 & 32.78 & 17.97 & 84.01\\
\hline
\end{tabular}
\end{center}
\end{table}

\begin{figure}[h]\label{fig:com_time}
\begin{center}
\includegraphics[width=0.75 \textwidth]{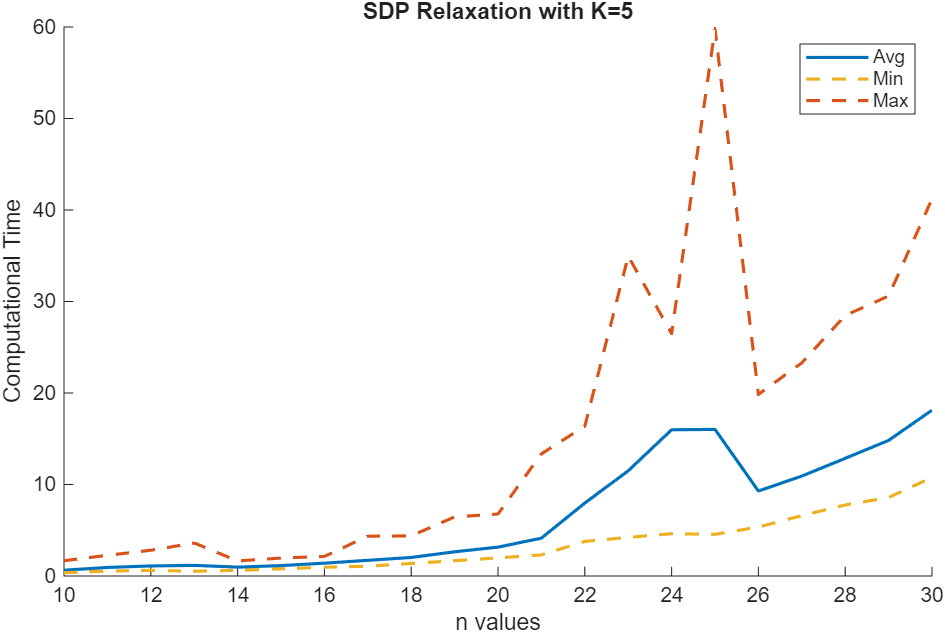}
% The figure caption is below the figure.
\end{center}
%\vspace{-1 in}
\caption{Computational Time of the SDP relaxation for $n\in\{10, 11, \dots, 30\}$ with $K=5$.}
\end{figure}

In the second experiment, we further investigate the computational efficiency of the SDP relaxation by considering a larger problem sizes $n=\{10, 11, \dots, 30\}$ with $K=5$. For each $n$, we generate 100 independent instances and record the computational time of the SDP relaxation. The results are presented in Figure \ref{fig:com_time}. As shown in the figure, the SDP relaxation has a substantially slower growth trend compared to the \eqref{EXP-LP} (whose growth trend can be seen in \Cref{tab:performance}. This demonstrates the efficiency of the SDP relaxation for larger problem instances.

\section{Conclusions}
\label{sec:conclusion}
In this work we have discussed completely positive formulations to compute the largest expected value of maximum of $K$ general quadratic functions over a moment ambiguity set. We look at three possible scenarios of support information. The problem is NP-hard for cases of support of random variables in $\{0,1\}^n $ and $\mathbb{R}_n^+$ while it is polynomial time solvable for support in $\mathbb{R}_n$. We propose a completely positive program for solving the two NP-hard cases and SDP for the polynomial time solvable case. Our numerical results demonstrate the optimality gap in the bound for the case of Bernoulli support where the exact formulation is computable using an exponential sized linear program. In future it is of interest to study the quality of bounds for specific structured forms of the quadratic functions. 

\begin{acknowledgements}
The research of the first and third authors are partially funded by MOE Tier 2 Grant MOET2EP20124 titled "Distributionally robust submodular optimization: Theory, algorithms and applications".

The research of the second author is partially funded by the Anusadhan National Research Foundation grants CRG/2022/001510 titled "Enhancing Robustness of Markov Decision Processes" and ANRF/ARGM/2025/001951/MTR titled ``Optimization with Submodular Functions under Uncertainty''. 
\end{acknowledgements}

\appendix  %This command ends the counting of sections.
\begin{appendix}
\section{The formulation of (EXP-LP) for the Bernoulli Support}
\begin{align*}\label{EXP-LP}
    \tag{EXP-LP} \quad \max& \sum_{k\in[K]}\sum_{i=1}^n\sum_{j=1}^n A_k(i,j)P(\tilde{\xi}_i=1,\tilde{\xi}_j=1, k^*(\bm{\tilde{\xi}})=k)\\
   &+\sum_{k\in[K]}\sum_{i=1}^n b_k(i)P(\tilde{\xi}_i=1, k^*(\bm{\tilde{\xi}})=k)+\sum_{k\in[K]}c_kP(k^*(\bm{\tilde{\xi}})=k)\\
   \text{s. t.}& \\
   & \underline{\gamma_i}\leq \sum_{k\in[K]}P(\tilde{\xi}_i=1, k^\ast(\bm{\tilde \xi})=k)\leq \bar{\gamma}_i \quad \forall i\in[n], \\
   & \sum_{k\in[K]}P(\tilde{\xi}_i=1, \tilde{\xi}_j=1, k^\ast(\bm{\tilde \xi})=k)\geq\Sigma_{ij} \quad \forall (i,j)\in\binom{[n]}{2}, \\
   & \sum_{k\in[K]}\sum_{\bm{\tilde{\xi}}\in\{0,1\}^n} P(\tilde{\xi}_1,\dots,\tilde{\xi}_n, k^\ast(\bm{\tilde \xi})=k)=1,\\
   & P(\tilde{\xi}_i=1, k^\ast(\bm{\tilde \xi})=k) =\sum_{\bm{d}\in\{0,1\}^{n-1}} P\left(\tilde{\xi}_i=1, \bigcap_{j\in[n]/\{i\}}(\tilde{\xi}_j=d_j), k^\ast(\bm{\tilde \xi})=k\right) \quad \forall i\in[n], \forall k\in[K],\\
   & P(\tilde{\xi}_i=1, \tilde{\xi}_j=1, k^\ast(\bm{\tilde \xi})=k) = \sum_{\bm{d}\in\{0,1\}^{n-2}} P\left(\tilde{\xi}_i=1, \tilde{\xi}_j=1, \bigcap_{l\in[n]/\{i,j\}}(\tilde{\xi}_l=d_l),k^\ast(\bm{\tilde \xi})=k\right) \quad \forall (i,j)\in\binom{[n]}{2}, \forall k\in[K],\\
   & P(k^\ast(\bm{\tilde \xi})=k) = \sum_{\bm{d}\in\{0,1\}^n} P(\tilde{\xi}_1=d_1,\dots, \tilde{\xi}_n=d_n, k^\ast(\bm{\tilde \xi})=k) \quad \forall k\in[K],\\
   & P(\tilde{\xi}_1, \dots, \tilde{\xi}_n, k^\ast(\bm{\tilde \xi}))\geq 0 \quad \forall \bm{\tilde{\xi}}\in\{0,1\}^n, \forall k\in[K].
\end{align*}

\end{appendix}
%References
\bibliographystyle{plain}
\bibliography{references}

\end{document}